\documentclass[12pt]{amsart}
\usepackage{amssymb,amsmath,euscript}
\usepackage[hypertex]{hyperref}
\usepackage{color}
\catcode`\@=11
\def\@evenfoot{\rule{0pt}{20pt}[\today] \hfill}
\def\@oddfoot{\rule{0pt}{20pt}\hfill [\today]}

\swapnumbers
\newtheorem{theorem}{Theorem}[section]
\newtheorem{corollary}[theorem]{Corollary}

\newtheorem{lemma}[theorem]{Lemma}
\newtheorem{proposition}[theorem]{Proposition}

\newtheorem*{theoremA}{Theorem~A}
\newtheorem*{theoremB}{Theorem~B}
\newtheorem*{theoremC}{Theorem~C}
\newtheorem*{theoremD}{Theorem~D}
\newtheorem*{theoremE}{Theorem~E}
\newtheorem*{theoremF}{Theorem~F}

\theoremstyle{definition}

\newtheorem{odstavec}[theorem]{}
\newtheorem{example}[theorem]{Example}
\newtheorem{remark}[theorem]{Remark}
\newtheorem{definition}[theorem]{Definition}

\newtheorem*{formality-theorem}{Formality}

\DeclareMathOperator{\id}{id} \def\RRR{\mathbb R} \def\frakS{\mathfrak S}
\DeclareMathOperator{\sign}{sign} \def\deviace{{\EuScript D}}
\newcommand{\oset}[2]{\overset{#1}{#2}}

\newcommand{\kriv}[3]{R({#1},{#2})({#3})}
 \def\bbbR{{\mathbb R}} \def\bbR{\bbbR}
\def\Tr{{\it Tr\/}}  \def\Rn{{\bbR^n}} 
\def\sqdots{. \hskip -.01em .\hskip -.01em .}
\def\deltah{\delta_{\rm h}}\def\Graphname{\Lambda}
\def\vforder{\hbox{${\it ord}_{\rm vf}$}}

\def\deltav{\delta_{\rm v}}\def\Tot{{\rm Tot\/}} \def\Zh{Z_{\rm h}}
\def\cyclsum#1{{\raisebox{-.1em}{\Large $\circ$}\hskip-1.2em\sum_#1}}
\def\cyclsumvf#1{{\raisebox{-.1em}{\Large $\circ$}\hskip-1.5em\sum_#1}}
\def\textcyclsum{{\raisebox{-.05em}{\large $\circ$}\hskip-.85em\sum}}
\def\bbB{{\mathbb B}} \def\Ker{{\it Ker\/}} \def\Im{{\it Im\/}} 
\def\iden{{\it id\/}} \def\gO{{\EuScript O\/}} \def\vt{\vartheta}
 \def\Span{{\it Span\/}}
\def\dr#1#2{\frac{\partial #1}{\partial #2}} 
\def\Grd #1#2{\Gr^{#1}[#2](d)} \def\Kr{{\EuScript K}}
\def\Gr{{\EuScript {G}\rm r}} \def\Con{{\it Con\/}}\def\Vert{{\it Vert\/}}

\def\Cont{{\EuScript {C}\rm ont}}

\def\Rada#1#2#3{#1_{#2},\dots,#1_{#3}} \def\lot{{\it l.o.t.}}
\def\bbox{{\raisebox {.1em}{\rule{.6em}{.6em}} \hskip .1em}} 
\def\sbbox{{\raisebox {.1em}{\rule{.4em}{.4em}} \hskip .1em}} 
\def\pa{\partial}  \def\Con{{\it Con\/}}
\def\Nat{{\mathfrak {Nat}\/}} 
\def\rada#1#2{#1,\ldots,#2}
\def\prodrada#1#2{#1 \cdots #2}
\def\Rada#1#2#3{#1_{#2},\dots,#1_{#3}}

\def\cases#1#2#3#4{
                  \left\{
                         \begin{array}{ll}
                           #1,\ &\mbox{#2}
                           \\
                           #3,\ &\mbox{#4}
                          \end{array}
                   \right.
}

\def\anchor{\unitlength .25cm
\begin{picture}(1,1.4)(-1,-.7)
\put(-.45,.55){\makebox(0,0)[cc]{$\sbbox$}}
\put(-.5,-.8){\vector(0,1){1.2}}
\end{picture}}

\def\anteanchor#1{\unitlength .25cm
\begin{picture}(1.5,1)(-1,-.7)
\put(-.45,-.55){\makebox(0,0)[cc]{$\bullet$}}
\put(-.45,-.55){\makebox(0,0)[lt]{\hskip .2em \scriptsize $#1$}}
\put(-.5,-.2){\vector(0,1){1.2}}
\end{picture}}

\title[COMBINATORIAL DIFFERENTIAL GEOMETRY]
{COMBINATORIAL DIFFERENTIAL GEOMETRY AND IDEAL BIANCHI--RICCI IDENTITIES}

\author[Jany\v ska, Markl]{J.~Jany\v ska, M.~Markl}

\thanks{The first author was supported
by the Ministry of Education of the Czech Republic under the Project
MSM0021622409. The second author was supported 
   by the grant GA \v CR 201/08/0397 and by
   the Academy of Sciences of the Czech Republic,
   Institutional Research Plan No.~AV0Z10190503}

\address{
{\ }
\newline
Department of Mathematics and Statistics, Masaryk University
\newline
Kotl\'a\v{r}sk\'a 2, 611 37 Brno, The Czech Republic
\newline
E-mail: {\tt janyska@math.muni.cz}
\newline
{\ }
\newline
Mathematical Institute of the Academy 
\newline
{\v Z}itn{\'a} 25,
         115 67 Prague 1, The Czech Republic
\newline
E-mail: {\tt markl@math.cas.cz}
}

\keywords{Natural operator, linear connection, reduction theorem, graph.
}

\subjclass[2000]{20G05, 53C05, 58A32}

\begin{document}

\bibliographystyle{plain}
\baselineskip18pt plus 1pt minus 1pt
\parskip3pt plus 1pt minus .5pt

\maketitle

\begin{abstract} 
We apply the graph complex method of~\cite{markl:na} to vector fields
depending naturally on a set of vector fields and a linear symmetric
connection.  We characterize all possible systems of generators for
such vector-field valued operators including the classical ones given
by normal tensors and covariant derivatives. We also describe the size
of the space of such operators and prove the existence of an `ideal'
basis consisting of operators with given leading terms which satisfy
the (generalized) Bianchi--Ricci identities without the correction
terms.
\end{abstract}

\noindent 
{\bf Plan of the paper.} In Sections~\ref{sec1} and~\ref{sec2} we
recall classical reduction theorems and the Bianchi-Ricci
identities. The main results of this paper, Theorems~A--F, are formulated in
Section~\ref{sec3}. Sections~\ref{sec4},~\ref{sec5} and~\ref{sec6}
contain necessary notions and results of the graph complex theory and related
homological algebra. Section~\ref{sec7} provides proofs of the
statements of Section~\ref{sec3}.

\section{Classical reduction theorems}
\label{sec1}

In this paper, $M$ will always denote a smooth manifold. The letters
$X$, $Y$, $Z$, $U$, $V$,..., with or without indices, will denote
(smooth) vector fields on $M$. The local coefficients of a vector
field $X$ are smooth functions $X^\lambda $ in coordinates $x^\lambda
$, $1 \leq \lambda \leq \dim(M)$, such that $ X = X^\lambda
\frac{\pa}{\pa x^{\lambda }} $ where, as usual, the summation over
repeated indices is assumed.  We also consider a linear connection $\Gamma
$ on $M$ with Christoffel symbols $\Gamma^\lambda _{\mu \nu }$, $1
\leq \lambda ,\mu ,\nu \leq \dim(M)$, see, for example,~\cite[Section
III.7]{kobayashi-nomizu}. The letter $R$ will denote the curvature
$(1,3)$-tensor field of $\Gamma $, the symbol $\nabla$ the
covariant derivative with respect to $\Gamma $, and by $\nabla^{(r)}$
we will denote the sequence of iterated covariant derivatives up to
order $r$, i.e.  $\nabla^{(r)} = (\id,\nabla,\dots,\nabla^r)$.  In
this paper we assume the connection 
$\Gamma$ to be {\em symmetric\/} (also called
{\em torsion-free\/}), i.e.\ $\Gamma^\lambda _{\mu \nu}=\Gamma^\lambda
_{\nu \mu}$.  The case of non-symmetric connections will be addressed
in a forthcoming paper.

It is well-known that natural (that is invariant with respect to chosen local 
coordinates, i.e.~coordinate-independent) 
operators of linear symmetric
connections on manifolds and of tensor fields
which have values in tensor fields
can be factorized through the curvature tensors, their covariant derivatives, 
given tensor fields and their covariant derivatives. These results are known as
the first (the operators of connections only) and the second reduction
theorems.  

\smallskip
Let us quote the original Schouten's formulation of the {\em first
reduction theorem\/}.

\begin{theorem}\label{Th1.1}
{\em \cite[p.\ 164]{Sch54}} 
{All differential concomitants of a symmetric
connexion are ordinary concomitants
of $R^{.}_\nu{}^._\mu{}^._\lambda{}^\kappa$
and its covariant derivatives.}
\end{theorem}

Similarly we have the {\em second reduction  theorem}.

\begin{theorem}\label{Th1.2}
{\em \cite[p.\ 165]{Sch54}} 
{All differential concomitants of a set of quantities
$\Phi_1,\dots, \Phi_d$ (indices suppressed) and the symmetric
connexion $\Gamma^\kappa_{\mu\lambda}$ are ordinary concomitants
of ${\Phi}_1, \dots,{\Phi}_d$, $R^{.}_\nu{}^._\mu{}^._\lambda{}^\kappa$
and their covariant derivatives.}
\end{theorem}

Let us recall that a {\em differential concomitant\/} is a polynomial
coordinate-independent operator with values in tensor fields depending on
certain order derivatives of input fields, while  {\em ordinary
concomitants\/} are zero order operators obtained from input fields by
tensorial operations, i.e.~by tensor products, permutations of indices
and contractions.\label{concom}

Proofs of the above reduction theorems use  normal coordinates
of $\Gamma$ centered at a point $x_0$ of $M$, see~\cite{VebTho23}. 
In such coordinates,
\begin{equation}
\label{0}
\Gamma^\lambda_{\mu\nu}(x^\rho) = x^\rho N{}^\lambda_{\rho\mu\nu}(x_0) +
\frac1{2!}\, x^{\rho_1}x^{\rho_2} N{}^\lambda_{\rho_1\rho_2\mu\nu}(x_0)
+ \cdots,
\end{equation}
where $N_n := (N{}^\lambda_{\rho_1\cdots\rho_{n-2}\mu\nu})$, $n\ge 3$, are
the {\it normal tensors} satisfying the following identities:
\begin{equation}
\label{1}
N_n(X_{\sigma(1)},\ldots,X_{\sigma(n-2)},  X_{n-1}, X_n ) - N_n(\Rada X1n) = 0
\end{equation}
for any permutation $\sigma$ of $(n-2)$ indices,
\begin{equation}
\label{2}
N_n(X_{1},\dots,X_{n-2},  X_{n-1}, X_n ) - 
N_n(X_{1},\dots,X_{n-2},  X_{n}, X_{n-1} ) = 0
\end{equation}
and 
\begin{equation}
\label{3}
\sum_{\sigma\in \Sigma_n} N_n(X_{\sigma(1)},\dots,X_{\sigma(n)}) =0,
\end{equation}
where $\Sigma_n$ denotes the permutation group of $n$ elements.  The
independence of a differential concomitant on given local coordinates
implies that any differential concomitant of $\Gamma$ is an
ordinary concomitant of the normal tensors $N_n$, $n\geq 3$. This
result is known as the {\it replacement theorem}, see~\cite{ThoMic27}.
The first reduction theorem now follows from the fact that each $N_n$ can
be expressed as a linear combination, with real coefficients, of the
covariant derivatives of order $(n-3)$ 
of the curvature tensor $R$ of
$\Gamma$ and a tensor field constructed from covariant derivatives of
orders $\le (n-4)$ of $R$, \cite[p. 162]{Sch54}, i.e., if we denote by
$$
(\nabla^i R)(X_1,\dots,X_i)(X_{i+1},X_{i+2})(X_{i+3})
$$
the $(1,i+3)$ tensor field $\nabla^i R$ evaluated on vector fields
$(X_1,\dots,X_{i+3})$, one can write
\begin{equation*}
N_n(X_1,\dots,X_n) = \sum_{\sigma\in\Sigma_n} A_\sigma
(\nabla^{n-3}R)(X_{\sigma(1)},\dots,X_{\sigma(n-3)})
(X_{\sigma(n-2)},X_{\sigma(n-1)})(X_{\sigma(n)})
+ l.o.t.\,,
\end{equation*}
with some $A_\sigma \in \mathbb{R}$, where $l.o.t.$ is a $(1,n)$-type ordinary
concomitant constructed from $\nabla^{(n-4)}R$.

The second reduction theorem can be proved similarly.

In the language of natural bundles and
natural operators, differential concomitants are polynomial natural
differential operators, see, for instance,
\cite{kolar-michor-slovak,KruJan90, Nij72, terng:AMJ78}.    
It is proved in
\cite[Section 28]{kolar-michor-slovak} 
that the above reduction
theorems are true for all natural differential operators, not only for
polynomial ones.

\smallskip

We will study polynomial natural differential operators
on vector fields and symmetric linear connections with values in vector fields.
By the second reduction theorem such operators of order $r$ with respect to
the vector fields are given by their $r$-th order covariant derivatives.
So, the minimal order with respect to the connection is $(r-1)$ but this order
can be higher. We have:

\begin{corollary}\label{Cr1.3}
All polynomial vector fields depending naturally on vector fields
${X}_1, \dots,{X}_d$ (in order $r$) and a symmetric linear connection
$\Gamma $ (in order $s\ge (r-1)$) are obtained by tensorial
constructions from the covariant derivatives (up to the order $(s-1)$)
of the curvature tensor of $\Gamma$, vector fields ${X}_1,
\dots,{X}_d$ and their covariant derivatives up to the order $r$, i.e.
\[
X(\Gamma;{X}_1, \dots,{X}_d) = 
	X(\nabla^{(s-1)}R;\nabla^{(r)}{X}_1,\dots,\nabla^{(r)}{X}_d)\,.
\]
\end{corollary}

\noindent 
{\bf Results of the paper.}  Classical reduction theorems, as
Theorems~\ref{Th1.1} and~\ref{Th1.2} above, describe systems of
operators which generate all operators of a given type. One usually
also gives a list of relations between these generators. For example,
the generating system $\oset{S}{\nabla}{}^{(s-1)} R$ for operators of
a symmetric connection given by the symmetrized covariant derivatives
of the curvature tensor satisfies the classical (higher order) Bianchi
identities~\eqref{Eq2.15} with non-vanishing right-hand sides.

In this paper we characterize all generating systems of natural
operators from a set of vector fields and a connection with values in
vector fields, see Theorems~A, B and C of
Section~\ref{sec3}. Theorem~D of the same section is a uniqueness
result for presentations in a given generating system. 
Theorem~E then states that, for each
choice of the leading terms of the generating operators, there exists
a streamlined, `ideal' version of these operators satisfying the (generalized)
Bianchi-Ricci identities without the right-hand sides. 
Finally, in Theorem~F, we calculate the size of
the spaces of natural operators studied in the paper.

\vskip .5em

\noindent 
{\bf Methods of the paper.}  Classical proofs of reduction theorems,
as those given in~\cite[Section III.7]{Sch54}, as well as proofs that
use techniques of natural operators, see~\cite[Section
28]{kolar-michor-slovak}, are based on technically complicated
calculations in coordinates.  The proofs given in this paper combine
the classical methods of normal coordinates with the graph complex
method proposed by the second author in~\cite{markl:na} which is
independent on local calculations and analysis.

While the `classical' methods are suited for proving that a certain
system of operators generate all operators of a given type, the
graph-complex method is particularly useful for analyzing the
uniqueness of expressing an operator via the generating ones. 
Therefore, the strength of the graph complex method will be particularly
manifest in Theorems~B,~D and~F of Section~\ref{sec3}.

Let us close this section by recalling 
that the graph complex method is a sophisticated version
of the `abstract tensor calculus.' It represents geometric
quantities, such as coordinates of a tensor field, via vertices of
graphs, with graphs playing the role of contraction schemes for composed
geometric objects. The coordinate independence of such expressions is
characterized by the vanishing of a graph differential. This 
brings our method into the realm of homological algebra. 
See~\cite{markl:na,markl:ig} for details.

\section{Classical  Bianchi and Ricci identities}
\label{sec2}
\setcounter{equation}{0}

If $Q$ is a tensor field of the type $(1,k)$, $k\ge 0$, 
we denote by 
$(\nabla^rQ)(X_1,\dots,X_r,Z_1,\dots,Z_k)$ 
the value of the $r$th covariant derivative of $Q$
evaluated on $(r+k)$ vector fields $X_j,Z_i$, $j=1,\dots,r$,
$i=1,\dots,k$, i.e.
$$
(\nabla^rQ)(X_1,\dots,X_r,Z_1,\dots,Z_k) = 
X^{\nu_1}\cdots X^{\nu_r}Z^{\mu_1}\cdots Z^{\mu_k}\, 
\nabla_{\nu_1}\cdots \nabla_{\nu_r}
Q^\lambda_{\mu_1\dots \mu_k}\frac{\partial}{\partial x^\lambda }\,.
$$
Let us denote by
$$ (\oset{S}{\nabla}{}^r Q)(X_1,\dots,X_r,Z_1,\dots,Z_k) =
\frac1{r!}\sum_{\sigma \in \Sigma_r} ({\nabla}^r
Q)(X_{\sigma(1)},\dots,X_{\sigma(r)},Z_1,\dots,Z_k)\,
$$
and
$$ (\oset{A}{\nabla}{}^r Q)(X_1,\dots,X_r,Z_1,\dots,Z_k) =
\frac1{r!}\sum_{\sigma \in \Sigma_r} (-1)^{\sign(\sigma)} ({\nabla}^k
Q)(X_{\sigma(1)},\dots,X_{\sigma(r)},Z_1,\dots,Z_k)\,
$$
the symmetrized and the antisymmetrized 
$r$th covariant derivatives of $Q$, respectively. 

Then we have the Ricci
identity
\begin{align}\label{Ricci1}
(\oset{A}{\nabla}{}^2Q)(Y,X, & Z_1,\dots,Z_k) 
	= 
- \frac12\left[\rule{0pt}{1em}\right.  R(X,Y)(Q(Z_1,\dots,Z_k) )
\\ \nonumber
&
- Q(R(X,Y)(Z_1),Z_2,\dots,Z_k) 
 - \dots -
	Q(Z_1,\dots,R(X,Y)(Z_k))\left.\rule{0pt}{1em}\right]\,.
\end{align}
{}From the Ricci identity \eqref{Ricci1} we obtain
\begin{align}\label{Eq2.2}
(\nabla^r Q)(X_1,\dots,X_r,Z_1,\dots,Z_k) & = (\oset{S}{\nabla}{}^r Q)(X_1,\dots,X_r,Z_1,\dots,Z_k)
\\ \nonumber
&\quad
+ pol(X_1,\dots,X_r,Z_1,\dots,Z_k)\,,
\end{align}
where $pol$ is a  $(1,r+k)$-type ordinary concomitant constructed from $\nabla^{(r-2)}Q$ and 
$\nabla^{(r-2)}R$. 
  
For vector fields we have
\begin{align*}
	(\nabla Y)(X) 
& = 
	\nabla_{X} Y\,,\qquad
	(\nabla^2 Z)(Y,X) 
 = 
	\nabla_{Y} (\nabla_{X} Z) - \nabla_{\nabla_{Y} X} Z\,.  
\end{align*}
Identity \eqref{Eq2.2} now has the form
\begin{align}
\label{Eq2.4}
(\nabla^r Z)(X_1,\dots,X_r) & = (\oset{S}{\nabla}{}^r Z)(X_1,\dots,X_r)
+ pol(X_1,\dots,X_r)\,,
\end{align}
where $pol$ is a  $(1,r)$-type ordinary concomitant constructed from $\nabla^{(r-2)}Z$ and 
$\nabla^{(r-2)}R$. 

For the curvature tensor we have the antisymmetry identity
\begin{equation}
\kriv{X}{Y}{Z} = - \kriv{Y}{X}{Z}\,,
\end{equation}
the 1st Bianchi identity
\begin{equation}\label{Eq2.6}
\cyclsumvf {{X,Y,Z}} \kriv{X}{Y}{Z} = 0\,,
\end{equation}
and the 2nd Bianchi identity
\begin{equation}\label{Eq2.7}
\cyclsumvf {{U,X,Y}} (\nabla R)(U)(X,Y)(Z) = 0\,,
\end{equation}
where $\textcyclsum$ denotes the cyclic summation over the 
indicated vector fields.
Identity \eqref{Eq2.2} for $R$ has the form
\begin{align}\label{Eq2.11}
(\nabla{}^r R)(U_1,\dots,U_r)(X,Y)(Z) 
&  = 
(\oset{S}\nabla{}^r R)(U_1,\dots,U_r)(X,Y)(Z) 
\\ \quad\nonumber
& + pol(U_1,\dots,U_r,X,Y,Z)\,, 
\end{align}
where $pol$ is a $(1,r+3)$-type ordinary concomitant constructed from 
$\nabla^{(r-2)}R$. 

\begin{remark}\label{Rm2.15}
The $r$th order, $r\ge 2$, covariant derivative
$\nabla^rR$ satisfies the identities obtained by the covariant
derivatives of the 1st and the 2nd Bianchi identities,~i.e.
\begin{equation}
\cyclsumvf {{X,Y,Z}} (\nabla^r R)(U_1,\dots,U_r)(X,Y)(Z) = 0\,,\qquad 
\cyclsumvf {{U_r,X,Y\hbox{\hglue-0.6ex}}} (\nabla^r R)(U_1,\dots,U_r)(X,Y)(Z) = 0\,.
\end{equation}
The symmetrized  $r$th order covariant derivative
$\oset{S}{\nabla}{}^r R$ however
satisfies only the identity obtained by the covariant derivative of the 1st
Bianchi identity. 

So, for the
symmetrized higher order covariant derivatives of the curvature tensor, we have the following 
(higher order) antisymmetry identity
\begin{equation}
 (\oset{S}\nabla{}^r R)(U_1,\dots,U_r)(X,Y)(Z) + (\oset{S}\nabla{}^r R)(U_1,\dots,U_r)(Y,X)(Z) = 0\,, 
\end{equation}
the (higher order) classical 1st Bianchi identity 
\begin{equation}
\cyclsumvf {{X,Y,Z}} (\oset{S}\nabla{}^r R)(U_1,\dots,U_r)(X,Y)(Z) = 0\,
\end{equation}
and the (higher order) classical 2nd Bianchi identity with a non-vanishing right hand side
\begin{equation}\label{Eq2.15}
\cyclsumvf {{U_r,X,Y\hbox{\hglue-0.6ex}}} (\oset{S}\nabla{}^r R)(U_1,\dots,U_r)(X,Y)(Z) = 
- \cyclsumvf {{U_r,X,Y\hbox{\hglue-0.6ex}}} pol(U_1,\dots,U_r,X,Y,Z)\,, 
\end{equation}
where $pol$ is a $(1,r+3)$-type ordinary concomitant
 from \eqref{Eq2.11}.
\end {remark}

\section{Main results}
\label{sec3}
\setcounter{equation}{0}

\begin{odstavec}
{\bf Operators we consider.}
\label{types}
Let $\Con$ be the natural bundle functor of torsion-free linear
connections~\cite[Section~17.7]{kolar-michor-slovak} and $T$ the tangent
bundle functor. We will consider natural differential operators $\gO :
\Con \times T^{\otimes d} \to T$ acting on a linear connection and $d$
vector fileds, $d \geq 0$, which are linear in the vector fields
variables, and which have values in vector fields. We will denote the space of
natural operators of this type by $\Nat(\Con \times T^{\otimes d},T)$.
Some typical operators $\Con \times T^{\otimes d} \to T$
are recalled in Example~\ref{Lie} on page~\pageref{Lie}.

Define the {\em vf-order\/} (vector-field order)\label{huhu}
resp.\ the {\em c-order\/} (connection order)\label{haha}
of a differential operator $\gO : \Con
\times T^{\otimes d} \to T$ as the order of $\gO$ 
in the vector field variables, resp.\ the connection variable.
\end{odstavec}

\begin{odstavec}{\bf Traces.} 
\label{traces}
Let $\gO$ be an operator acting on vector fields $X_1,\ldots,X_d$ and
a connection $\Gamma$, with values in vector fields.  Suppose that
$\gO$ is a linear order $0$ differential operator in $X_i$ for some $1
\leq i \leq d$. This means that the local formula $O(\Gamma,\Rada
X1d)$ for $\gO$ is a linear function of the coordinates of $X_i$ and
does not contain derivatives of the coordinates of $X_i$. In this
situation we define $\Tr_i(\gO) \in \Nat(\Con \times T^{\otimes (d-1)},R)$ 
as the operator with values in the bundle $R$ of
smooth functions given by the local formula
\begin{eqnarray*}
\lefteqn{
\Tr_i(O)(\Gamma,\Rada X1{i-1},\Rada X{i+1}d) :=} \hskip 1em
\\ 
&& \hskip 1em
\mbox{Trace}(O (\Gamma,\Rada X1{i-1},-,\Rada X{i+1}d) : \Rn \to \Rn).
\end{eqnarray*}
It is easy to see that $\Tr_i(\gO)$  is well defined. Whenever we
write $\Tr_i(\gO)$ we tacitly assume that the trace makes sense, 
i.e.~that $\gO$ is linear order $0$ differential
operator in $X_i$. 
\end{odstavec}

\begin{odstavec}
{\bf Compositions} 
\label{comp}
Let $\gO' : \Con \times T^{\otimes d'} \to T$ and $\gO'' : \Con \times
T^{\otimes d''} \to T$ be operators as in~\ref{types}. Assume that
$\gO'$ is a linear order $0$ differential operator in $X_i$ for some
$1 \leq i \leq d'$. In this situation we define the composition $\gO'
\circ_i \gO'' : \Con \times T^{\otimes (d'+ d'' -1)} \to T$ as the
operator obtained by substituting the value of the operator $\gO''$
for the vector-field variable $X_i$ of $\gO'$. As in~\ref{traces}, by
writing $\gO' \circ_i \gO''$ we signalize that $\gO'$ is of order $0$
in $X_i$.
\end{odstavec}

\begin{odstavec}
\label{iter}
Throughout this section, by an {\em iteration\/} of differential
operators we understand applying a finite number of the following
`elementary' operations:
\begin{itemize}
\item[(i)] permuting the vector-fields inputs of a differential operator $\gO$,

\item[(ii)] taking the pointwise linear
combination $k' \cdot \gO' + k'' \cdot \gO''$, $k',k''\in \mathbb{R}$,

\item[(iii)] performing the substitution $\gO' \circ_i \gO''$, and

\item[(iv)] taking the pointwise product $\Tr_{i}(\gO') \cdot \gO''$.
\end{itemize}

We of course assume that the operations in (ii) and (iii) make sense,
see~\ref{traces} and~\ref{comp}. There are `obvious' relations between
the above operations. The operations $\circ_i$ in (iii) satisfy the `operadic'
associativity and are compatible with permutations in~(i),
see properties (1.9) and (1.10) 
in~\cite[Definition~II.1.6]{markl-shnider-stasheff:book}. Other
`obvious' relations are the commutativity of the trace, $\Tr_j(\gO'
\circ_i \gO'') = \Tr_i(\gO''\circ_j \gO')$ and its `obvious'
compatibility with permutations of~(i).

The iteration defined above provides a coordinate-independent definition 
of an ordinary concomitant recalled on page~\pageref{concom}, 
i.e.~an operator $\gO$ is an iteration of
operators $\Rada \gO1N$ if and only if it is an ordinary concomitant
of $\Rada \gO1N$.   
\end{odstavec}

\begin{odstavec}
\label{shuffles}
Let us consider, for $n \geq 3$, the induced representation 
$E^0(n) := {\rm Ind}^{\Sigma_n}_{\Sigma_{n-2}\times \Sigma_2} 
({{\mathbf 1}}_{n-2}\times {{\mathbf 1}}_2)$,
where ${{\mathbf 1}}_{n-2}$ (resp.~${{\mathbf 1}}_{2}$)
is the trivial representation of the
symmetric group $\Sigma_{n-2}$ (resp.~$\Sigma_2$). Elements of
$E^0(n)$ are linear combinations
\begin{equation}
\label{konec}
\sum_{\sigma \in {\rm Ush}(n-2,2)} \alpha_\sigma 
\cdot (1_{n-2}\times 1_2) \sigma,
\end{equation}
where $1_{n-2}\times 1_2 \in {{\mathbf 1}}_{n-2}\times {{\mathbf
1}}_2$ is the generator, $\alpha_\sigma \in {\mathbb R}$, and $\sigma$
runs over all $(n-2,2)$-unshuffles which are, by definition,
permutations $\sigma \in \Sigma_n$ such that $\sigma(1) < \cdots
<\sigma(n-2)$ and $\sigma(n-1) < \sigma(n)$. Let $E^1(n)$ be the
trivial $\Sigma_n$-module ${\mathbf 1}_n$ and $\vartheta_E : E^0(n) \to E^1(n)$
the equivariant map that sends the generator $1_{n-2} \times 1_2 \in {{\mathbf
1}}_{n-2}\times {{\mathbf 1}}_2$ to $-1_n \in {{\mathbf 1}}_n$. The
reason for this notation and sign convention will became clear in
Section~\ref{cohomology}.

Define finally $\Kr(n) \subset E^0(n)$ to be
the kernel of the map $\vartheta_E$.  It is clear that  $\Kr(n)$
consists of all
expressions~(\ref{konec}) such that
\begin{equation}
\label{Konec}
\sum_{\sigma \in {\rm Ush}(n-2,2)} \alpha_\sigma  = 0.
\end{equation}
\end{odstavec}

\begin{theoremA}
\label{A}
Let $D_n(\Gamma,\Rada X1n)$, $n \geq 3$, be differential operators in 
$\Nat(\Con \times T^{\otimes n},T)$ whose local expressions are
\begin{equation}
\label{oolong}
D_n^\omega\left(\Gamma^\lambda_{\mu\nu},
\rada{X^{\delta_1}_1}{X^{\delta_n}_n}\right) = \sum_{\sigma \in 
{\rm Ush}(n-2,2)} \alpha_\sigma \cdot 
X^{\rho_1}_{\sigma(1)}\cdots X^{\rho_n}_{\sigma(n)} 
\frac{\pa^{n-2}\Gamma^\omega_{\rho_{n-1}\rho_n}}
{\prodrada{\pa x^{\rho_1}}{\pa x^{\rho_{n-2}}}} + \lot
\end{equation}
where $\lot$ is an expression of differential order $< n-2$, and
$\{\alpha_\sigma\}_{\sigma \in \Sigma_n}$ are real constants such that the
element~(\ref{konec}) generates the $\Sigma_n$-module $\Kr(n)$ introduced
in~\ref{shuffles} (which in particular means that~(\ref{Konec}) is
satisfied).

Let also $V_n(\Gamma,\Rada X1n)$,
$n \geq 1$, be  differential operators in 
$\Nat(\Con \times T^{\otimes n},T)$ of the form
\[
V_n^\omega\left(\Gamma^\lambda_{\mu\nu},
\rada{X^{\delta_1}_1}{X^{\delta_n}_n}\right) =
X^{\rho_1}_1\cdots X^{\rho_{n-1}}_{n-1}  
\frac{\pa^{n-1} X^{\omega_n}}
{\pa x^{\rho_1}\cdots{\pa x^{\rho_{n-1}}}
}
 + \lot,
\]
where $\lot$ is an expression of differential order 
$< n-1$.  

Suppose, moreover, that the operator $D_n(\Gamma,\Rada X1n)$ is of
vf-order~$0$ and the operator $V_n(\Gamma,\Rada X1n)$ of order $0$ in
$\Rada X1{n-1}$.  Then each differential operator $\gO : \Con \times
T^{\otimes d} \to T$ is an iteration, in the sense of~\ref{iter}, of
the operators $\{D_n\}_{n \geq 3}$ and $\{V_n\}_{n \geq 1}$.
\end{theoremA}

Theorem~A, as well as other statements in this Section, are proved in
Section~\ref{proofs}. Observe that necessarily $V_1(\Gamma,X) = X$, so
we may safely discard $V_1$ from the list of `generating' operators
and consider $V_n$'s only for $n \geq 2$.

\begin{remark}
\label{rem1}
It is a simple exercise on the Littlewood-Richardson rule
that the $\Sigma_n$-module $\Kr(n)$ of~\ref{shuffles} decomposes as 
$\Kr(n) \cong \bigoplus_{\lambda} V_\lambda$,
with the summation taken over all two-column Young diagrams
$\lambda = (\lambda_1,\lambda_2)$ with $\lambda_1 \geq 2$, $0 <
\lambda_2 \leq \lambda_1$, such that $\lambda_1 + \lambda_2 = n$, and where
$V_\lambda$ is the irreducible representation corresponding to
$\lambda$. Since all irreducible factors of $\Kr(n)$ have
multiplicity one, an element $x \in \Kr(n)$ is a $\Sigma_n$-generator
if and only if $\pi_\lambda(x) \not=0$ for each
projection $\pi_\lambda :  \Kr(n)\to V_\lambda$. Therefore the assumption
of Theorem~A on the coefficients $\alpha_\sigma$ 
can in principle be checked.

The operator $D_n(\Gamma,\Rada X1n)$ in Theorem A has vf-order 0 so it can be
interpreted as a $(1,n)$-tensor field $D_n(\Gamma)$ naturally
depending on $\Gamma$ (with c-order $(n-2)$). Similarly $V_n(\Gamma,\Rada X1n)$
can be considered as a $(1,n-1)$-tensor field $V_n(\Gamma,X_n)$
naturally given by $\Gamma$ and $X_n$ (with order $(n-1)$ with respect
to $X_n$).

Then the set $\{D_n(\Gamma)\}_{n \geq 3}$ and $\{V_n(\Gamma,X_i)\}_{n
\geq 1}$, $i = 1,\dots, d$, is a new system of generating operators
for natural vector fields from Corollary \ref{Cr1.3}.
\end{remark}

There are two `preferred' choices of the leading terms of the
operators $D_n$ in Theorem~A, the expression
\begin{equation}
\label{K}
X^{\rho_1}_1\cdots X^{\rho_n}_n \frac{\pa^{n-3}}
{\prodrada{\pa x^{\rho_1}}{\pa x^{\rho_{n-3}}}}
\left(
\frac{\pa \Gamma^\omega_{\rho_{n-2}\rho_n}}{\pa x^{\rho_{n-1}}}
- \frac{\pa  \Gamma^\omega_{\rho_{n-1}\rho_n}}{\pa x^{\rho_{n-2}}}
\right)
\end{equation}
and the expression
\begin{equation}
\label{N}
\left( 
\sum_{\sigma \in {\rm Ush}(n-2,2)} \frac{n(n-1)}{2}
X^{\rho_1}_{\sigma(1)}\cdots X^{\rho_n}_{\sigma(n)} 
-  X^{\rho_1}_1\cdots X^{\rho_n}_n
\right)
\frac{\pa^{n-2}\Gamma^\omega_{\rho_{n-1}\rho_n}}
{\prodrada{\pa x^{\rho_1}}{\pa x^{\rho_{n-2}}}}.
\end{equation}

The leading term~(\ref{K}) is given by the following choice of the
coefficients in~(\ref{oolong}):
\[
\alpha_\sigma :=
\left\{\begin{array}{rl}
-1& \mbox{if $\sigma$ is the identity},
\\
1& \mbox{if $\sigma$ is the 2-cycle that interchanges $(n-2)$ and
  $(n-1)$, and}
\\
0& \mbox {otherwise.} 
\end{array}
\right.
\]
It is fairly easy to prove that the corresponding element
in~(\ref{konec}) generates $\Kr(n)$. The element~(\ref{konec})
corresponding to~(\ref{N}) can be written as
\[
(1_{n-2}\times 1_2) 
\left(\sum_{\sigma \in {\rm Ush}(n-2,2)} \frac{n(n-1)}{2}
\sigma
- {\rm id}
\right),
\] 
so it is the image of the generator $1_{n-2}\times 1_2$ of $E^0(n)$
under the projection $E^0(n) \twoheadrightarrow \Kr(n)$. This
immediately implies that it generates $\Kr(n)$.

Expression~(\ref{K}) is the leading term of the iterated covariant
derivative of the curvature resp.~of its
streamlined version whose existence is proved in  
Theorem~E. The leading
term of~(\ref{N}) is that of the ``normal tensors'' of~(\ref{0}), see
also Example~\ref{e88} below.

\begin{example}
\label{e88}
Operators having the form required by Theorem~A exist. One may, for
instance, take
\[
K_n(\Gamma,\Rada X1n) := 
(\nabla^{n-3} R)({X_1}, \cdots ,X_{n-3})(X_{n-2},X_{n-1})(X_n),\ n \geq 3,
\]
for the operators $D_n$ and
\[
V_n(\Gamma,\Rada X1n) := 
(\oset{S}{\nabla}{}^{n-1} X_n)(X_{1}, \cdots ,X_{n-1}),\ n \geq 2.
\]
It is obvious that the leading term of $K_n$ is
expression~(\ref{K}). Another realization of the operators $D_n$ is
provided by the normal tensors $N_n$ recalled in~(\ref{0}) whose leading
term is~(\ref{N}).
\end{example}

For an operator $\gO \in \Nat(\Con \times T^{\otimes n},T)$ and a
permutation $\sigma \in \Sigma_n$ we denote by $\gO \sigma$ the
operator obtained by permuting the vector-field variables $\Rada X1n$
of $\gO$ according to $\Sigma_n$. This action extends, by linearity, into
a right action of the group ring $\RRR [\Sigma_n]$. We will denote
$\gO c$ the result of the action of $c \in \RRR [\Sigma_n]$ on $\gO$. 
The following theorem characterizes all possible systems
of generating operators.

\begin{theoremB}
Assume that $\dim(M) \geq 3$.  Let $U_n(\Gamma,\Rada X1n)
\in \Nat(\Con \times T^{\otimes n},T)$, $n \geq 3$, be operators of
vf-order $0$ and of c-order $(n-2)$. Then the following two conditions
are equivalent.

(i) Each operator $\gO \in \Nat(\Con \times T^{\otimes
d},T)$ of vf-order $0$ and c-order $(n-2)$ is an iteration of the
operators $\{U_u\}_{u \leq n}$.  

(ii) For each $n \geq 3$, there are elements $c,\Rada c1n$ of the
group ring $\RRR [\Sigma_n]$ 
such that the leading term of the operator
\begin{equation}
\label{j12}
D_n := U_n c + \sum_{1 \leq j \leq n}\Tr_j(U_n c_j) X_j \in \Nat(\Con
\times T^{\otimes n},T)
\end{equation}
is of the form required by Theorem~A.
\end{theoremB}

\begin{example}
Consider the operator 
\[
U_3(\Gamma,X,Y,Z)  := R(X,Y)(Z) + \Tr(R(-,Z)(X)Y +  \Tr(R(-,Z)(Y) X.
\]
Then clearly $\frac 12\{U(\Gamma,X,Y,Z) - U(\Gamma,Y,X,Z)\}$ equals the
curvature $R(X,Y)(Z)$, so $D_3$ defined by~(\ref{j12}) with $c
= \frac 12 \tau_{12}$, $c_1 = c_2 = c_3 = 0$, where $\tau_{12}$ is the
permutation $(1,2,3) \mapsto (2,1,3)$, has the leading term
required by Theorem~A. By Theorem~B, $U_3$ defined above can be a
member of a generating series of operators.
\end{example}

\begin{example}
Let us illustrate the necessity of the assumption $\dim(M)
\geq 3$ in Theorem~B.
Let $U_3 \in  \Nat(\Con \times T^{\otimes 3},T)$ be the operator
defined by
\[
U_3(\Gamma,X,Y,Z) := X \cdot \Tr(R(Y,-)(Z)) + Y \cdot \Tr(R(-,X)(Z)).
\]
The leading term of this operator equals
\[
X^\omega Y^\mu Z^\nu \frac{\pa \Gamma^\lambda_{\mu\nu}}{\pa x^\lambda}
-Y^\omega X^\mu Z^\nu \frac{\pa \Gamma^\lambda_{\mu\nu}}{\pa
  x^\lambda}
+ Y^\omega X^\lambda Z^\nu \frac{\pa \Gamma^\mu_{\mu \nu}}{\pa
  x^\lambda}
-X^\omega Y^\lambda Z^\nu \frac{\pa \Gamma^\mu_{\mu \nu}}{\pa x^\lambda},
\]
so it is clearly \underline{not} of the form required by Theorem~A.
On the other hand, it can be verified by a straightforward calculation
that on a $2$-dimensional manifold,
\[
U_3(\Gamma,X,Y,Z) =  R(X,Y)(Z) + \lot
\]
therefore, in dimension 2, the operator $U_3$ can be a part of a
generating series of operators.
\end{example}

\begin{theoremC}
Assume that $\dim(M) \geq 2d-1$ and that $\{D_n\}_{n \geq 3}$, 
$\{V_n\}_{n \geq 1}$ be as in Theorem~A.
Let $\gO : \Con \times T^{\otimes d} \to T$ be a differential operator
of the vf-order  $a \geq 0$. Then it has an iterative representation
with the following property. Suppose that an 
additive factor of this iterative 
representation of $\gO$ via $\{D_n\}_{n \geq 3}$ and
$\{V_n\}_{n \geq 2}$ contains
$\rada {V_{q_1}}{V_{q_t}}$, for some $\Rada q1t \geq 2$,  
$t \geq 0$. Then
\[
q_1+\cdots +q_t \leq a +t.
\]
In particular, if $\gO$ is of vf-order~$0$, then there exists an 
iterative representation that uses only $\{D_n\}_{n \geq 3}$.
\end{theoremC}

Notice that one can prove the particular case of Theorem~C for
operators of vf-order~$0$ without the  $\dim(M) \geq 2d-1$ assumption
by a simple modification of the `classical' proof of Theorem~A given
on page~\pageref{J}. We, however, do not know how to use the classical
reduction techniques to prove Theorem~C in full generality.

\begin{example}
\label{nedele}
It is clear that $[X,Y] = \nabla_X Y - \nabla_Y X$ i.e.~if $V_2$ is as
in Example~\ref{e88},
\[
[X,Y] = V_2(\Gamma,X,Y) - V_2(\Gamma,Y,X).
\]
This shows that the individual summands of an iterative representation
of an operator $\gO$ may depend on the connection though the operator
$\gO$ does not. This fact was used in~\cite{KruJan90}, in an other context,
as the method of an {\em auxiliary connection\/}. 
Later, it was proved in~\cite{markl:na} that operators in 
$\Nat(\Con \times T^{\otimes d},T)$ that do not depend on the
connection are iterations of the Lie bracket
of vector fields.
\end{example}

\begin{example}
\label{Jar}
We show that there is, in general, no relation
between the c-order of a
differential operator and the c-order of its iterative
representation. 
We have
\begin{eqnarray*}
[X,[Y,Z]] &=& V_3(\Gamma,X,Y,Z) - V_3(\Gamma,X,Z,Y) +
V_2(\Gamma,V_2(\Gamma,X,Y),Z) 
\\ 
&& - V_2(\Gamma,V_2(\Gamma,X,Z),Y) -
V_2(\Gamma,V_2(\Gamma,Y,Z),X) 
\\
&&+ V_2(\Gamma,V_2(\Gamma,Z,Y),X) + \lot
\end{eqnarray*}
While the c-order of $[X,[Y,Z]]$ is~$0$, the operators $V_3$ in
the right hand side are of the c-order~$1$.
\end{example}

For $n \geq 3$, $\sigma \in \Sigma_n$ and $D_n$ as in Theorem~A,
denote by $D_n \sigma$ the operator obtained from $D_n$ by permuting
the vector fields variables according to $\sigma$. This notation
clearly extends to the action of an element $\frakS$ of the group ring
$\RRR [\Sigma_n]$.

\begin{definition}
\label{debra}
We say that $\frakS \in \RRR [\Sigma_n]$ is a {\em quasi-symmetry\/} of
an operator $D_n$ in~(\ref{oolong}) if
\[
(\sum_{\sigma \in \Sigma_n} \alpha_\sigma \sigma )\frakS = 0 
\]
in the group ring $ \RRR [\Sigma_n]$.
We say that  $\frakS$ is a {\em symmetry\/} of $D_n$ if $D_n \frakS = 0$.
\end{definition}

A  quasi-symmetry  $\frakS$ of $D_n$, by definition, annihilates its 
leading term, therefore $D_n \frak S$ is an operator of
c-order $\leq (n-3)$ that does not use the derivatives of the vector
field variables. We can express this fact by writing
\begin{equation}
\label{dev}
D_n \frak S(\Gamma,X_1,\ldots,X_n) = 
\deviace^\frakS_n(\Gamma,X_1,\ldots,X_n)
\end{equation}
where $\deviace^\frakS_n \in \Nat(\Con \times T^{\otimes n},T)$
($\deviace$ abbreviating ``deviation'') is a degree $\leq (n-3)$
operator which is, by Theorem~C, an iteration of the operators $D_u$
with $3 \leq u \leq n-1$ (no $V_n$'s). By definition, $\frakS$ is a symmetry
of $D_n$ if and only if $\deviace^\frakS_n =0$.  We will see,
in~\ref{odst} below, that~(\ref{dev}) offers a conceptual explanation
of the Bianchi and Ricci identities.

A similar discussion can be made also for the operators 
$V_n$, $n \geq 1$. Since the leading term of $V_n$ is fully symmetric in
$\Rada X1{n-1}$,
\begin{equation}
\label{b5}
V_n(\Gamma,X_{\omega(1)},\ldots,X_{\omega(n-1)},X_n) -
V_n(\Gamma,\Rada X1n) = 
\deviace V^{\omega}_n(\nabla,\Rada X1n)
\end{equation}
for any $n \geq 2$ and $\omega \in \Sigma_{n-1}$, where $\deviace V^{\omega}
\in \Nat(\Con \times T^{\otimes n},T)$ is an order $\leq (n-2)$
differential operator.

The following theorem states that the iteration of Theorem~A is unique
up to identities~(\ref{dev}),~(\ref{b5}) and the `obvious' relations.

\begin{theoremD}
On manifolds of dimension $\geq 2d-1$, the iteration expressing an
operator in $\Nat(\Con \times T^{\otimes d},T)$ via
$\{D_n\}_{n \geq 3}$ and $\{V_n\}_{n \geq 2}$ is unique up to
relations~(\ref{dev}) with $\frakS$ running over all quasisymmetries of
$D_n$, relation~(\ref{b5}), and the `obvious' relations among
elementary operations.  
In particular,~(\ref{dev}) and~(\ref{b5}) are the only 
(quasi)symmetries of the operators 
$\{D_n\}_{n \geq 3}$ and $\{V_n\}_{n \geq 2}$.
\end{theoremD}

\begin{odstavec}\label{odst}
{\bf Bianchi and Ricci identities.}  
The leading term~(\ref{K}) enjoys the following symmetries:
\begin{itemize}
\item[(s1)] 
the antisymmetry in $X_{n-2}$ and $X_{n-1}$,
\item[(s2)]
the cyclic symmetry in $X_{n-2}$, $X_{n-1}$ and $X_n$, 
\item[(s3)]
for $n \geq 4$, the cyclic symmetry in $X_{n-3}$, $X_{n-2}$ and $X_{n-1}$, and
\item[(s3)]
for $n \geq 4$, the total symmetry in $X_1,\ldots,X_{n-3}$. 
\end{itemize}
We leave as an exercise to express these symmetries via appropriate
elements of the group ring $\RRR[\Sigma_n]$. 
It is not difficult to prove that (s1)--(s4) generate all symmetries
of~(\ref{K}).

Let $D_n$ be an operator of the form~(\ref{oolong}) with the leading
term~(\ref{K}). The (anti)symmetry (s1) leads to the equation
\begin{equation}
\label{b1}
D_n\left(\Gamma,X_1,\sqdots,X_{n-2},X_{n-1},X_n\right)+
D_n\left(\Gamma,X_1,\sqdots,X_{n-1},X_{n-2},X_n\right)
= \deviace^{\rm as}_n(\Gamma,X_1,\sqdots,X_n),
\end{equation}
where the natural differential operator $\deviace^{\rm as}_n$ 
of order $\leq (n-3)$  can be interpreted as the deviation
from antisymmetry of $D_n$. Similarly, (s2) leads to
\begin{equation}
\label{b2}
\cyclsum {\sigma}
D_n\left(
\Gamma,X_1,\ldots,X_{\sigma(n-2)},X_{\sigma(n-1)},X_{\sigma(n)}
\right)
= {\deviace'_n}^{\rm cycl}(\Gamma,\Rada X1n),
\end{equation}
where $\textcyclsum$ is the cyclic summation over the indicated
indices and ${\deviace'_n}^{\rm cycl} \in \Nat(\Con \times T^{\otimes n},T)$ is
an order $\leq (n-3)$ differential operator.  In the same manner, for $n
\geq 4$, (s3) gives
\begin{equation}
\label{b3}
\cyclsum {\sigma}
D_n\left(
\Gamma,X_1,\ldots,X_{n-4},X_{\sigma(n-3)},X_{\sigma(n-2)},X_{\sigma(n-1)},X_n
\right)= {\deviace''_n}^{\rm cycl}(\Gamma,\Rada X1n),
\end{equation}
for some order $\leq (n-3)$ operators 
${\deviace''_n}^{\rm cycl} \in \Nat(\Con \times T^{\otimes n},T)$.
Finally, for $n \geq 5$, the symmetry (s4) implies
that
\begin{eqnarray}
\label{b4}
\lefteqn{\hskip -14cm
D_n\left(
\Gamma,X_{\omega(1)},\ldots,X_{\omega(n-3)},X_{n-2},X_{n-1},X_n
\right)
-D_n\left(
\Gamma,X_1,\ldots,X_n
\right) =}
\\
\nonumber
\hskip 12cm = \deviace^{\omega}_n(\Gamma,\Rada X1n)
\end{eqnarray}
for each $n \geq 4$ and a permutation $\omega \in \Sigma_{n-3}$, with
some order  $\leq (n-3)$ operators 
$\deviace^{\omega}_n\in \Nat(\Con \times T^{\otimes n},T)$. 

Observe that, by Theorem~C, the right hand sides 
of~(\ref{b1})--(\ref{b4}) are iterations of the operators $D_u$
with $3 \leq u \leq n-1$ (no $V_n$'s).
We call~(\ref{b5})--(\ref{b4}) the {\em Bianchi-Ricci\/} identities for the
operators $\{D_n\}_{n \geq 3}$ and $\{V_n\}_{n \geq 2}$.
\end{odstavec}

\begin{example}\label{Ex3.8}
It is clear that not only the leading term of the operator $V_n$ of
Example~\ref{e88}, but the operator itself is fully symmetric in
$\Rada X1{n-1}$, therefore~(\ref{b5}) for this operator 
is satisfied with trivial right hand side, 
\[
V_n(\Gamma,X_{\omega(1)},\ldots,X_{\omega(n-1)},X_n) -
V_n(\Gamma,\Rada X1n)= 0
\]
for any $n \geq 2$ and $\omega \in \Sigma_{n-1}$. 

Let us inspect
symmetries fulfilled by the operators $K_n$ of Example~\ref{e88}.  

For $n=3$,~(\ref{b1}) reduces to the standard antisymmetry of the
curvature tensor,
\begin{equation}
\label{6}
R(X,Y)(Z) + R(Y,X)(Z)=0
\end{equation}
and, for $n\geq 4$,~(\ref{b1}) is the iterated covariant derivative
of~(\ref{6}). Therefore $K_n$ satisfies~(\ref{b1}) with the trivial
right hand side,
\[
K_n\left(\Gamma,X_1,\ldots,X_{n-2},X_{n-1},X_n\right)+
K_n\left(\Gamma,X_1,\ldots,X_{n-1},X_{n-2},X_n\right)=0.
\]
Similarly, for $n = 3$,~(\ref{b2}) means 
the vanishing of  the cyclic sum,
\begin{equation}
\label{7}
R(X,Y)(Z) + R(Y,Z)(X) + R(Z,X)(Y) = 0,
\end{equation}
which is the classical 1st Bianchi identity \eqref{Eq2.6} of a torsion-free
connection. For 
$n \geq 4$, (\ref{b2}) is the iterated covariant derivative of~(\ref{7}),
therefore
\[
\cyclsum {\sigma}
K_n\left(
\Gamma,X_1,\ldots,X_{n-3},X_{\sigma(n-2)},X_{\sigma(n-1)},X_{\sigma(n)}
\right)
= 0.
\]

For $n = 4$,
the left hand side of~(\ref{b3}) 
means the cyclic sum,
\begin{equation}
\label{7a}
(\nabla R)(U)(X,Y)(Z) + (\nabla R)(X)(Y,U)(Z) + (\nabla R)(Y)(U,X)(Z) ,
\end{equation}
and, by the classical 2nd Bianchi identity \eqref{Eq2.7}, it is
satisfied with the vanishing right hand side.  For $n \geq 5$,
(\ref{b3}) is the iterated covariant derivative of~(\ref{7a}),
therefore
\[
\cyclsum {\sigma}
K_n\left(
\Gamma,X_1,\ldots,X_{\sigma(n-4)},X_{\sigma(n-3)},X_{\sigma(n-2)},X_{n}
\right)
= 0.
\]

On the other hand, for $n \ge 5$, the left hand side of~(\ref{b4}) is  
given by the
covariant derivatives of the Ricci identities, and it is \underline{nonzero}.
\end{example}

\begin{odstavec}\label{normal}
{\bf Symmetries of the normal tensors.}
The leading term~(\ref{N}) of the normal tensor $N_n$ has the
following symmetries:
\begin{itemize}
\item[(s1)]
the full symmetry in $\Rada X1{n-2}$,
\item[(s2)]
the symmetry in $X_{n-1}$ and $X_n$, and
\item[(s3)]
the symmetry described by $\frakS := \sum_{\sigma \in \Sigma_n} \sigma$. 
\end{itemize}
One can prove that (s1)--(s3) generate
all symmetries of~(\ref{N}).
Equations~(\ref{1})--(\ref{3}) then say that these symmetries of the
leading term in fact extend to symmetries of the operator~$N_n$.

\end{odstavec}

The following theorem shows that for each choice 
of the leading terms there exist particularly nice operators $\{D_n\}_{n
\geq 3}$ of Theorem~A.

\begin{theoremE}
For each choice of the leading terms satisfying~(\ref{Konec}), 
there exist `ideal' operators $\{J_n\}_{n \geq 3}$ of the
form~(\ref{oolong}), for which all the
``generalized'' Bianchi-Ricci identities~(\ref{dev}) are satisfied without
the right hand sides. In other words, all quasisymmetries, in the
sense of Definition~\ref{debra}, 
are actual symmetries of the operators $\{J_n\}_{n \geq 3}$.
\end{theoremE}

Observe that, in Theorem~E, we do not assume that the
element~(\ref{konec}) related to the leading term generates $\Kr(n)$,
we only assume that it belongs to the kernel of the map $\vartheta_E$.

\begin{example}
\label{ex09} 
As we already saw in Example~\ref{Ex3.8}, the operators $V_n$
introduced in Example~\ref{e88} are `ideal' in the sense that all
their quasi-symmetries are also their symmetries. This is not true for
the iterated covariant derivatives $K_n$ of the curvature tensor
(see again Example~\ref{Ex3.8}), neither for their `naively' symmetrized
versions ${\widetilde K}_n :=  \oset{S}{\nabla}{}^{(n-3)} R$, $n \geq
3$. The ideal versions $J_n$ of these operators (both having the same leading
term) which exist by Theorem~E, can be constructed by modifying
${\widetilde K}_n$ as
\[
J_n = {\widetilde K}_n + P_n,
\]
where $P_n(\Gamma,\Rada X1n)$ is a c-order $\le (n-3)$ iteration of
the operators $\{{\widetilde K}_u\}_{3 \leq u \leq n-1}$. While
clearly $P_3= P_4 = 0$, i.e.
\[
J_3 = R\,,\quad J_4 = \nabla R,
\]
the calculation of the correction term $P_n$ is, for $n \geq 5$, a
nontrivial task. To give the reader the taste of the complexity of
the calculation, we write an explicit formula for $P_5$:
\begin{align*}
	P(\Gamma,& U,V,X,Y,Z) 
 =
- \tfrac12 \left\{2 R(U,R(X,Z)(Y))(V) - 2 R(U,R(Y,Z)(X))(V)
\rule{0pt}{1.2em}\right.
\\
& \quad	
	+ 2 R(V,R(X,Z)(Y))(U) 
- 2 R(V,R(Y,Z)(X))(U)
	+ 2 R(U,R(X,Y)(V))(Z) 
\\
& \quad	
+ 2 R(V,R(X,Y)(U))(Z)
	+ R(X,R(U,Z)(V))(Y) + R(X,R(V,Z)(U))(Y)
\\ 
& \quad	
	- R(Y,R(U,Z)(V))(X) - R(Y,R(V,Z)(U))(X)
	+ R(U,R(X,Z)(V))(Y) 
\\
& \quad	
+ R(V,R(X,Z)(U))(Y)
	- R(U,R(Y,Z)(V))(X) - R(V,R(Y,Z)(U))(X)
\\ 
& \quad
	+ R(Y,R(X,U)(V))(Z) + R(Y,R(X,V)(U))(Z)
	- R(X,R(Y,U)(V))(Z) 
\\
& \quad	
- R(X,R(Y,V)(U))(Z)
	+ R(Y,R(X,U)(Z))(V) + R(Y,R(X,V)(Z))(U)
\\ 
& \quad	
	- R(X,R(Y,U)(Z))(V) - R(X,R(Y,V)(Z))(U)
	+ R(X,R(U,Z)(Y))(V) 
\\
& \quad	\left.
+ R(X,R(V,Z)(Y))(U)
	- R(Y,R(U,Z)(X))(V) - R(Y,R(V,Z)(X))(U)
\rule{0pt}{1.2em}\right\}\,.
\end{align*}
\end{example}

\begin{remark}
Let as remark that the ideal basis $\{J_n\}_{n\ge 3}$ of the type 
discussed in Example~\ref{ex09} can be constructed using the
normal tensors $\{N_n\}_{n\ge 3}$ of~\eqref{0} as
\begin{align}
J_n(X_1,& \dots,X_{n-3},X_{n-2},X_{n-1},X_n) := 
\\ \nonumber
& = N_n(X_1,\dots,X_{n-3},X_{n-2},X_{n-1},X_n)
- N_n(X_1,\dots,X_{n-3},X_{n-1},X_{n-2},X_n)\,.
\end{align}
It is indeed easy to see that the operators $J_n$ defined in this way
have the same leading terms as the operators $K_n$ and that
identities~\eqref{1}--\eqref{3} imply the 
identities~\eqref{b1}--\eqref{b4} with trivial right-hand sides.

On the other hand, the normal tensor $N_n$ is, for each $n \geq 3$,  
a linear combination
\begin{equation}
N_n(X_1,\dots,X_n) = \sum_{\sigma\in \Sigma_n}  A_\sigma \,
J_{n}(X_{\sigma(1)},\dots,X_{\sigma(n)}) ,
\end{equation}
where the real coefficients $A_\sigma$ are determined by requiring that 
identities~\eqref{b1}--\eqref{b4} with trivial right-hand sides imply
identities~\eqref{1}--\eqref{3}.
\end{remark}

The size of the space of natural operators 
$\Con \times T^{\otimes d} \to T$ is described in the last theorem of
this section:

\begin{theoremF}
On manifolds of dimension $\geq 2d-1$, the vector space $\Nat(\Con
\times T^{\otimes d},T)$ is isomorphic to the graph space $\Gr[\Kr](d)$
introduced on page~\pageref{kriz} of Section~\ref{cohomology}.
\end{theoremF}

\begin{example}
The calculation of the dimension of $\Gr[\Kr](d)$ as of a
vector space spanned by graphs is a purely combinatorial problem. For
instance, for $d=1$ we get $\dim(\Gr[\Kr](d)) = 1$, with the
corresponding natural operator the identity $X \mapsto X$.

One also easily calculates that, on manifolds of dimension $\geq 3$,
$\dim(\Nat(\Con \times T^{\otimes 2},T)) = \dim(\Gr[\Kr](d)) = 4$. The
corresponding generating operators are 
\[
\nabla_XY,\ \nabla_YX,\ X
\cdot \Tr (\nabla_- Y)\ \mbox {and}\ Y \cdot \Tr (\nabla_- X).
\]
\end{example}

\begin{figure}
\setlength{\unitlength}{1mm}
\mbox{\begin{picture}(110,85)
\put(60,78){\makebox(0,0)[cc]{$\{N_n\}_{n\ge 3}$}}
\put(60,73){\makebox(0,0)[cc]{identities of normal tensors}}
\put(49,65){\vector(-1,-2){25}}
\put(71,65){\vector(1,-2){25}}
\put(60,40){\makebox(0,0)[cc]{general $\{D_n\}_{n\ge 3}$}}
\put(60,35){\makebox(0,0)[cc]{general B.-R. identities}}
\put(60,45){\vector(0,1){20}}
\put(42,30){\vector(-1,-1){15}}
\put(78,30){\vector(1,-1){15}}
\put(23,10){\makebox(0,0)[cc]{ideal $\{J_n\}_{n\ge 3}$}}
\put(23,5){\makebox(0,0)[cc]{ideal B.-R. identities}}
\put(40,10){\vector(1,0){40}}
\put(24,15){\vector(1,2){25}}
\put(97,10){\makebox(0,0)[cc]{$\{\overset{s}{\nabla}{}^{n-3} R\}_{n\ge 3}$}}
\put(97,5){\makebox(0,0)[cc]{classical B.-R. identities}}
\put(80,10){\vector(-1,0){40}}
\put(96,15){\vector(-1,2){25}}
\end{picture}}
\caption{\label{fig66}Classical bases of operators as specializations of
  the universal one.}
\end{figure}
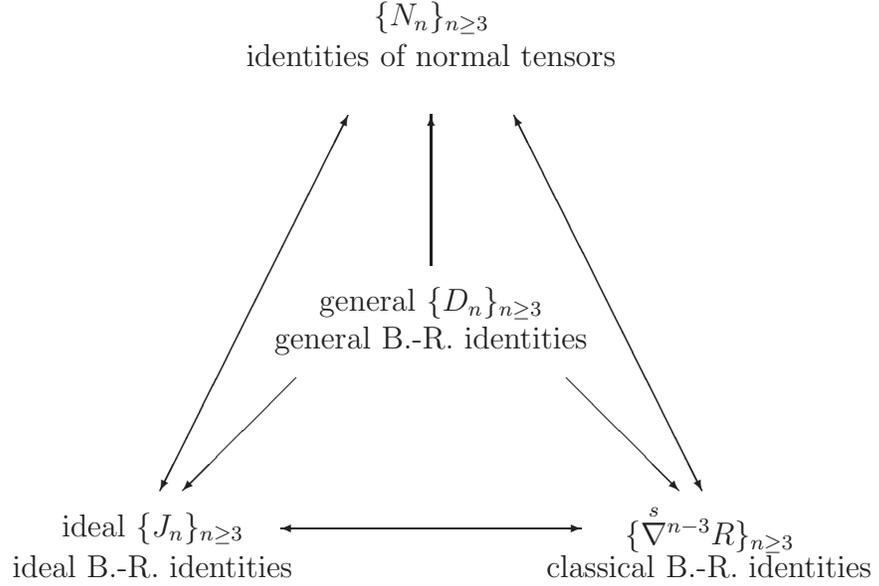

Results of this section characterize bases of natural operators in
$\Nat(\Con \times T^{\otimes d},T)$ and state some properties
of these bases. Various `classical' bases are then special cases of
these general bases. This is symbolically expressed by Figure~\ref{fig66}.

\section{Rules of the game}
\setcounter{equation}{0}
\label{rules}\label{sec4}

In this section whose bulk is taken from~\cite{markl:na} we recall the
graph complex describing natural differential operators from
$\Nat(\Con \times T^{\otimes d},T)$. The underlying graded vector
space of this complex is spanned by directed, not necessary
connected, graphs with three types of vertices plus one special vertex
called the anchor, see~\ref{o1}. The differential given by vertex
replacements is recalled in~\ref{o2}. Let us, however, start with an
example showing three typical operators from $\Nat(\Con \times
T^{\otimes d},T)$.

\begin{example}
\label{Lie}
The {\em Lie bracket\/}
$X,Y \mapsto [X,Y]$ is a natural operator that constructs from two
vector fields on $M$ a third one. For the purposes of this paper we
consider the Lie bracket as an operator in 
$\Nat(\Con \times T^{\otimes 2},T)$.
In local coordinates,
\begin{equation}
\label{e3334}
[X,Y] =
\left(X^\mu \frac{\pa Y^\lambda}{\pa x^\mu} - 
Y^\mu \frac{\pa X^\lambda}{\pa x^\mu}\right)\frac {\pa}{\pa x^\lambda}.
\end{equation}
The {\em covariant derivative\/}
$(\Gamma, X, Y) \mapsto \nabla_X Y$ is another
natural differential operator {}from $\Nat(\Con \times T^{\otimes 2},T)$. 
In local coordinates,
\begin{equation}
\label{e10}
\nabla_X Y = \left(\Gamma^\lambda_{\mu\nu} X^\mu Y^\nu 
+ X^\mu \frac{\pa Y^\lambda}{\pa x^\mu}
              \right) \dr{}{x^\lambda},
\end{equation}
where $\Gamma^i_{jk}$ are the Christoffel symbols.
The {\em curvature\/} $R \in \Nat(\Con \times T^{\otimes 3},T)$ of
$\Gamma$ is a~composition of the above operators,
\[
R(X,Y)(Z) := \nabla_{[X,Y]} Z - [\nabla_X,\nabla_Y] Z.
\]
\end{example}

\begin{odstavec}   
{\bf The graph complex.}
\label{o1}
In this paper, by a {\em graph\/} we mean a directed (i.e. with
oriented edges), not necessary connected, graph whose vertices are of
the types described below. Multiple edges, loops and wheels are
allowed.  Let us recall the construction of the graph complex
$\Gr^*(d)$ describing natural differential operators from $d$ vector
fields and a~torsion-free linear connection with values in vector
fields that are $d$-multilinear in the vector field variables.  Details
and motivations can be found in~\cite{markl:na} but observe that in
that paper we did not assume that the connection is torsion-free. The
degree $m$ part $\Gr^m(d)$ is spanned by graphs with precisely $d$
`black' vertices
\begin{equation}
\label{black}
\raisebox{-3em}{\rule{0pt}{0pt}}
b_u:= \hskip -6em 
\unitlength 4mm
\linethickness{0.4pt}
\begin{picture}(20,4)(10.5,18.8)
\put(20,20){\vector(0,1){2}}
\put(18,18){\vector(1,1){1.9}}
\put(22,18){\vector(-1,1){1.9}}
\put(19,18){\vector(1,2){.935}}
\put(20,20){\makebox(0,0)[cc]{\Large$\bullet$}}
\put(25,20){\makebox(0,0)[cc]{,\ $u \geq 0$,}}
\put(20,18){\makebox(0,0)[cc]{\scriptsize (\hskip 18mm)}}
\put(20.5,18){\makebox(0,0)[cc]{$\ldots$}}
\put(20,17){\makebox(0,0)[cc]{%
   $\underbrace{\rule{16mm}{0mm}}_{\mbox{\scriptsize $u$ inputs}}$}}
\end{picture}
\end{equation}
labelled by $X,Y,\ldots$, $\Rada X1d$ or $1,\ldots,d$,  some number of vertices
\begin{equation}
\label{nabla}
\unitlength 5mm
\begin{picture}(.7,2)(-.5,-.5)
\put(-.5,0){\makebox(0,0)[cc]{\Large$\nabla$}}
\put(2.5,0){\makebox(0,0)[cc]{, $u \geq 0$.}}
\put(-.5,.4){\vector(0,1){1}}
\put(.1,-1.25){\vector(-1,2){.5}}
\put(.59,-1.22){\vector(-1,1){.97}}
\unitlength 3mm
\put(-21,-20){
\put(19.4,18){\vector(1,4){0.475}}
\put(18.9,18){\vector(1,2){.95}}
\put(16,18){\vector(2,1){3.8}}
\put(17.75,18){\makebox(0,0)[cc]{$\ldots$}}
\put(17.7,16.6){\makebox(0,0)[cc]{%
   $\underbrace{\rule{10mm}{0mm}}_{\mbox{\scriptsize $u$ inputs}}$}}
\put(17.7,17.8){\makebox(0,0)[cc]{\scriptsize $( \hskip 12mm )$}}
\put(21.5,17.8){\makebox(0,0)[cc]{\scriptsize $( \hskip 4mm )$}}
}
\end{picture}
\raisebox{-2.5em}{\rule{0pt}{0pt}}
\end{equation}
precisely $m$ `white' vertices
\begin{equation}
\label{white}
\raisebox{-3em}{\rule{0pt}{0pt}}
\unitlength 4mm
\linethickness{0.4pt}
\begin{picture}(20,4)(10.5,18.8)
\put(20,20.3){\vector(0,1){1.7}}
\put(18,18){\vector(1,1){1.8}}
\put(22,18){\vector(-1,1){1.8}}
\put(19,18){\vector(1,2){.91}}
\put(20,20){\makebox(0,0)[cc]{\Large$\circ$}}
\put(25,20){\makebox(0,0)[cc]{,\ $u \geq 2$,}}
\put(20,18){\makebox(0,0)[cc]{\scriptsize (\hskip 18mm)}}
\put(20.5,18){\makebox(0,0)[cc]{$\ldots$}}
\put(20,17){\makebox(0,0)[cc]{%
   $\underbrace{\rule{16mm}{0mm}}_{\mbox{\scriptsize $u$ inputs}}$}}
\end{picture}
\end{equation}
and one vertex \anchor\ (the anchor). The braces $( \hskip 1em )$  in
the above pictures mean that the inputs they encompass are fully
symmetric, but we will usually omit these braces in the forthcoming text.

In the above graph complex, black vertices~(\ref{black}) correspond to
derivatives of vector field coordinates,
\[
X^\lambda_{(\mu_1,\ldots,\mu_u)}
:= \frac {\pa^u X^\lambda}{\pa x^{\mu_1} \cdots \pa
x^{\mu_u}},
\]
$\nabla$-vertices~(\ref{nabla}) to the derivatives of the
Christoffel symbols,
\[
{}_{(\omega_1,\ldots,\omega_u)}\Gamma^\lambda_{\mu\nu}:= 
\frac {\pa^u \Gamma^\lambda_{\mu\nu}}{\pa x^{\omega_1} \cdots \pa
x^{\omega_u}},
\] 
white vertices~(\ref{white}) correspond to generators of infinitesimal
symmetries and the anchor \anchor\ to the vector-field value of the
operator. 
\end{odstavec}

\begin{odstavec}
{\bf General connections and vector fields.} A simple modification of
material in~\ref{o1} describes operators from general, not necessarily
torsion-free, linear connections and vector fields into vector
fields. The corresponding graph complex has the same
vertices~(\ref{black}),~(\ref{white}) and \anchor\ as the graph
complex $\Gr^*(d)$, above but we do not assume that the
$\nabla$-vertex~(\ref{nabla}) is symmetric in the two rightmost
inputs.  The replacement rules are the same as for $\Gr^*(d)$ and the
obvious analog of Theorem~\ref{t2} holds.
\end{odstavec}

\begin{odstavec}
{\bf Anchored versus rooted graphs.} 
\label{ss2}
Before we proceed to the differential in the graph complex, we need to
make a couple of observations on the structure of our graphs.
All graphs we have been working with so far had
an anchor~\anchor.
For an arbitrary $k$, $0 \leq k \leq d$, 
denote by $\Gr^*(d)_k$  the subspace of  $\Gr^*(d)$
spanned by graphs 
$\Upsilon$ with a distinguished subset  
$\{\anteanchor1, \ldots, \anteanchor k\}$ 
of the set of black vertices~(\ref{black}) with $u=0$. 
Schematically such an $\Upsilon$ looks as
\begin{equation}
\label{Upsilon}
{
\unitlength=.800000pt
\thicklines
\begin{picture}(60.00,50.00)(0.00,20.00)
\put(30.00,30.00){\makebox(0.00,0.00){$\Upsilon$}}
\put(60.00,25.00){\makebox(0.00,0.00){.}}
\put(60.00,0.00){\makebox(0.00,0.00){$\bullet$}}
\put(10.00,0.00){\makebox(0.00,0.00){$\bullet$}}
\put(0.00,0.00){\makebox(0.00,0.00){$\bullet$}}
\put(34.00,0.00){\makebox(0.00,0.00){$\cdots$}}
\put(60.00,-7.00){\makebox(0.00,0.00)[t]{\scriptsize $k$}}
\put(10.00,-7.00){\makebox(0.00,0.00)[t]{\scriptsize $2$}}
\put(0.00,-7.00){\makebox(0.00,0.00)[t]{\scriptsize $1$}}
\put(31.00,70.00){\makebox(0.00,0.00){$\sbbox$}}
\put(60.00,0.00){\vector(-1,1){16.00}}
\put(10.00,0.00){\vector(2,3){9.00}}
\put(0.00,0.00){\vector(1,1){16.00}}
\put(30.00,50.00){\vector(0,1){20.00}}
\put(30.00,30.00){\circle{40.00}}
\end{picture}}
\raisebox{-2.2em}{\rule{0pt}{10pt}}
\end{equation}
An obvious right $\Sigma_k$-action on the space $\Gr^*(d)_k$ permutes
the labels of the distinguished vertices. Let $\Gr^*_k := \bigoplus_{d
\geq k}\Gr^*(d)_k$.  We will call graphs as in~(\ref{Upsilon}) {\em
anchored $k$-graphs\/}.

For each graph $\Upsilon \in \Gr^*_k$, there is a graph $\Upsilon^a$
with $k$ input edges and one output edge, obtained by amputating
$\sbbox$ from the anchor and $\bullet$'s from the input edges of the
distinguished vertices. For instance, the graph $\Upsilon$
in~(\ref{Upsilon}) gives the amputated graph
\begin{equation}
\label{Upsilona}
{
\unitlength=.800000pt
\thicklines
\begin{picture}(60.00,50.00)(0.00,20.00)
\put(25.00,25.00){\makebox(0.00,0.00)[lb]{$\Upsilon^a$}}
\put(34.00,0.00){\makebox(0.00,0.00){$\cdots$}}
\put(60.00,-7.00){\makebox(0.00,0.00)[t]{\scriptsize $k$}}
\put(10.00,-7.00){\makebox(0.00,0.00)[t]{\scriptsize $2$}}
\put(0.00,-7.00){\makebox(0.00,0.00)[t]{\scriptsize $1$}}
\put(60.00,0.00){\vector(-1,1){16.00}}
\put(10.00,0.00){\vector(2,3){9.00}}
\put(0.00,0.00){\vector(1,1){16.00}}
\put(30.00,50.00){\vector(0,1){20.00}}
\put(30.00,30.00){\circle{40.00}}
\end{picture}}
\raisebox{-2.2em}{\rule{0pt}{10pt}}
\end{equation}
with one output and $k$ numbered input edges. As in operad theory,
we will call the output edge of $\Upsilon^a$ the {\em root\/} and its
inputs the {\em legs\/}. Graphs as in~(\ref{Upsilona}) will be then
called {\em rooted $k$-graphs\/}. The operation $\Upsilon \mapsto
\Upsilon^a$ is clearly one-to-one and, when convenient, 
we will make no distinction
between anchored graphs and the associated rooted graphs.

Let $v$ be a vertex of a  graph 
$\Graphname  \in \Gr^* := \bigoplus_{d \geq 1} \Gr^*(d)$. 
Suppose that $v$ has $k$ input edges and let $\Upsilon
\in \Gr^*_k$ be a rooted $k$-graph enjoying the same $\Sigma_k$-symmetry as the
inputs of $v$.  One then may replace the
vertex $v$ by $\Upsilon$, by grafting the root of $\Upsilon$ to the
output edge of $v$ and the legs of $\Upsilon$ to the input edges of
$v$. We denote the result by  $\Graphname \circ_v
\Upsilon \in \Gr^*$ and call this operation the 
{\em vertex insertion\/} or {\em vertex replacement\/}. One can
clearly extend this by linearity to define $\Graphname \circ_v x$ for $x
\in \Gr^*_k$ a linear combination of graphs with appropriate symmetry.
\end{odstavec}

\begin{odstavec}
{\bf The differential.}
\label{o2}
Let us recall that the graph differential $\delta : \Gr^*(d) \to
\Gr^{*+1}(d)$ is the linear map acting on a graph $\Graphname \in
\Gr^m(d)$ by the formula
\begin{equation}
\label{e68}
\delta(\Graphname) = \sum_{v \in \Vert(\Graphname)}\epsilon_v \cdot
\Graphname \circ_v \delta(v) \in \Gr^{m+1}(d),\ m \geq 0,
\end{equation}
in which $\Vert(\Graphname)$ denotes the set of vertices of $\Graphname$ 
and $\Graphname \circ_v \delta(v)$
the result of replacing  the vertex $v$ by the sum $\delta(v)$ of rooted
graphs defined by 
\[
\raisebox{-4.2em}{\rule{0pt}{0pt}}
\hskip -2cm
\unitlength 3mm
\linethickness{0.4pt}
\begin{picture}(20,4)(10.5,18.8)
\put(16.8,19.5){\makebox(0,0){$\delta\left(\rule{0pt}{20pt}\right.$}}
\put(22.8,19.5){\makebox(0,0){$\left.\rule{0pt}{20pt}\right)$}}
\put(20,20.4){\vector(0,1){1.7}}
\put(18,18){\vector(1,1){1.8}}
\put(22,18){\vector(-1,1){1.8}}
\put(19,18){\vector(1,2){.91}}
\put(20,20){\makebox(0,0)[cc]{\Large$\circ$}}
\put(20.5,18){\makebox(0,0)[cc]{$\ldots$}}
\put(20,16.7){\makebox(0,0)[cc]{%
   $\underbrace{\rule{12mm}{0mm}}_{\mbox{\scriptsize $k$ inputs}}$}}
\end{picture}
\hskip -2.1cm :=  \sum_{s + u = k}
\hskip -1.2cm
\unitlength3mm
\linethickness{0.4pt}
\begin{picture}(20,4)(10.5,18.8)
\put(20,20.4){\vector(0,1){1.7}}
\put(18.1,18.1){\vector(1,1){1.7}}
\put(22,18){\vector(-1,1){1.8}}
\put(19,18){\vector(1,2){.91}}
\put(20,20){\makebox(0,0)[cc]{\Large$\circ$}}
\put(20.5,18){\makebox(0,0)[cc]{$\ldots$}}
\put(20.5,16.9){\makebox(0,0)[cc]{%
   $\underbrace{\rule{9mm}{0mm}}_{\mbox{\scriptsize $s$}}$}}
\end{picture}
\put(-22.5,-2.2){
\begin{picture}(20,4)(10.5,18.8)
\put(18,18){\vector(1,1){1.8}}
\put(22,18){\vector(-1,1){1.8}}
\put(19,18){\vector(1,2){.91}}
\put(20,20){\makebox(0,0)[cc]{\Large$\circ$}}
\put(20.5,18){\makebox(0,0)[cc]{$\ldots$}}
\put(20,17){\makebox(0,0)[cc]{%
   $\underbrace{\rule{12mm}{0mm}}_{\mbox{\scriptsize $u$}}$}}
\put(21.9,19){\makebox(0,0)[cc]{$\left(\rule{0pt}{20pt} \hskip 2.2cm
    \right)_{\rm ush}$}}
\end{picture}
}\hskip -4.5em  ,\ k \geq 2,
\]
for white vertices,
\begin{equation}
\label{e10bis}
\raisebox{-3.2em}{\rule{0pt}{0pt}}
\unitlength 3mm
\hskip -1cm
\linethickness{0.4pt}
\delta\left(
\begin{picture}(5,3)(17.5,18.8)
\put(20,20.3){\vector(0,1){1.7}}
\put(18,18){\vector(1,1){1.8}}
\put(22,18){\vector(-1,1){1.8}}
\put(19,18){\vector(1,2){.91}}
\put(20,20){\makebox(0,0)[cc]{\Large$\bullet$}}
\put(19.5,20.5){\makebox(0,0)[rb]{\scriptsize $X$}}
\put(20.5,18){\makebox(0,0)[cc]{$\ldots$}}
\put(20,16.7){\makebox(0,0)[cc]{%
   $\underbrace{\rule{12mm}{0mm}}_{\mbox{\scriptsize $k$ inputs}}$}}
\end{picture}
\right)
:= \sum_{s+u=k}
\hskip .7cm
\unitlength 3mm
\linethickness{0.4pt}
\begin{picture}(15,3)(15.5,18.8)
\put(20,20.3){\vector(0,1){1.7}}
\put(18,18){\vector(1,1){1.8}}
\put(22,18){\vector(-1,1){1.8}}
\put(19.1,18){\vector(1,2){.91}}
\put(20,20){\makebox(0,0)[cc]{\Large$\circ$}}
\put(20.5,18){\makebox(0,0)[cc]{$\ldots$}}
\put(20.5,16.8){\makebox(0,0)[cc]{%
   $\underbrace{\rule{9mm}{0mm}}_{\mbox{\scriptsize $s$}}$}}
\end{picture}
\put(-22.5,-2.2){
\begin{picture}(10,4)(10.5,18.8)
\put(18,18){\vector(1,1){1.8}}
\put(22,18){\vector(-1,1){1.8}}
\put(19.1,18){\vector(1,2){.91}}
\put(20,20){\makebox(0,0)[cc]{\Large$\bullet$}}
\put(19.5,20.5){\makebox(0,0)[rb]{\scriptsize $X$}}
\put(20.5,18){\makebox(0,0)[cc]{$\ldots$}}
\put(20,16.8){\makebox(0,0)[cc]{%
   $\underbrace{\rule{11mm}{0mm}}_{\mbox{\scriptsize $u$}}$}}
\put(21.7,19){\makebox(0,0)[cc]{$\left(\rule{0pt}{17pt} \hskip 2.1cm
\right)_{\rm ush}$}}
\end{picture}
}
\hskip -2cm
-
\hskip -1cm
\unitlength 3mm
\linethickness{0.4pt}
\begin{picture}(20,3)(10.5,18.8)
\put(20,20.3){\vector(0,1){1.7}}
\put(18.1,18.1){\vector(1,1){1.7}}
\put(22,18){\vector(-1,1){1.8}}
\put(19.1,18){\vector(1,2){.91}}
\put(20,20){\makebox(0,0)[cc]{\Large$\bullet$}}
\put(19.5,20.5){\makebox(0,0)[rb]{\scriptsize $X$}}
\put(20.5,18){\makebox(0,0)[cc]{$\ldots$}}
\put(20.5,16.8){\makebox(0,0)[cc]{%
   $\underbrace{\rule{9mm}{0mm}}_{\mbox{\scriptsize $s$}}$}}
\end{picture}
\put(-22.5,-2.2){
\begin{picture}(20,3)(10.5,18.8)
\put(18,18){\vector(1,1){1.8}}
\put(22,18){\vector(-1,1){1.8}}
\put(19.1,18){\vector(1,2){.91}}
\put(20,20){\makebox(0,0)[cc]{\Large$\circ$}}
\put(20.5,18){\makebox(0,0)[cc]{$\ldots$}}
\put(20,16.7){\makebox(0,0)[cc]{%
   $\underbrace{\rule{12mm}{0mm}}_{\mbox{\scriptsize $u$}}$}}
\put(21.7,19){\makebox(0,0)[cc]{$\left(\rule{0pt}{20pt} \hskip 2.1cm 
\right)_{\rm ush}$}}
\end{picture}
} \hskip -4.5em  ,\ k \geq 0,
\end{equation}
for black vertices and $\delta(\anchor) = 0$ for the anchor.  The
braces $(\hskip 1em )_{\rm ush}$ in the right hand sides
indicate that the summations over all $(u,s-1)$-unshuffles of
the inputs have been performed, see formulas~(29) and~(30)
of~\cite{markl:na} for details.
For convenience, we write explicitly formula~(\ref{e10bis}) for $k = 0,1,2$: 
\[
\raisebox{.4em}{\rule{0pt}{10pt}}
\raisebox{-.3em}{\rule{0pt}{10pt}}
\delta \left(\hskip 1em
\unitlength .4cm
\begin{picture}(0,0)(0,.2)
\put(-.5,0){\makebox(0,0)[cc]{\Large$\bullet$}}
\put(-.5,0){\vector(0,1){1.2}}
\end{picture} 
\right) =  0,\
\delta \left(\hskip 1em
\unitlength .4cm
\begin{picture}(0,0)(0,-.3)
\put(-.5,0){\makebox(0,0)[cc]{\Large$\bullet$}}
\put(-.5,0){\vector(0,1){1.2}}
\put(-.5,-1){\vector(0,1){.9}}
\end{picture} 
\right) =  
\unitlength 2.2 mm
\begin{picture}(4,1)(-2,-.6)
\put(0,0.6){\vector(0,1){1.7}}
\put(0,-2){\vector(0,1){1.78}}
\put(1.9,-1.9){\vector(-1,1){1.62}}
\put(1.7,-1.7){\makebox(0,0)[cc]{\Large$\bullet$}}
\put(0,0){\makebox(0,0)[cc]{\Large$\circ$}}
\end{picture},\
\delta \left(
\begin{picture}(4,1)(-2,-.6)
\put(0,0.3){\vector(0,1){1.7}}
\put(-1,-2){\vector(1,2){.91}}
\put(1,-2){\vector(-1,2){.91}}
\put(0,0){\makebox(0,0)[cc]{\Large$\bullet$}}
\end{picture}
\right) = -
\unitlength 2.2mm
\begin{picture}(4,1)(-2,0)
\put(0,0.6){\vector(0,1){2.5}}
\put(0,0.4){\vector(0,1){1}}
\put(-1,-2.3){\vector(1,2){.89}}
\put(1,-2.3){\vector(-1,2){.89}}
\put(0,1.5){\makebox(0,0)[cc]{\Large$\bullet$}}
\put(0,-0.2){\makebox(0,0)[cc]{\Large$\circ$}}
\end{picture} 
+
\unitlength 2.2mm
\begin{picture}(4,1)(-2,-.6)
\put(0,0.6){\vector(0,1){1.7}}
\put(-1.3,-2.7){\vector(1,2){1.2}}
\put(1.3,-2.7){\vector(-1,2){1.2}}
\put(1.3,-2.7){\vector(-1,2){.4}}
\put(.75,-1.7){\makebox(0,0)[cc]{\Large$\bullet$}}
\put(0,0){\makebox(0,0)[cc]{\Large$\circ$}}
\end{picture}
+
\unitlength 2.2mm
\begin{picture}(4,1)(-2,-.6)
\put(0,0.6){\vector(0,1){1.7}}
\put(-1,-2){\vector(1,2){.88}}
\put(1,-2){\vector(-1,2){.88}}
\put(2,-1){\vector(-2,1){1.7}}
\put(2,-1){\makebox(0,0)[cc]{\Large$\bullet$}}
\put(0,0){\makebox(0,0)[cc]{\Large$\circ$}}
\end{picture} \hskip .7em ,\ \ldots
\]

The replacement rule for the
$\nabla$-vertices is of the form
\begin{equation}
\label{e11}
\hskip 2cm
\unitlength .5cm
\begin{picture}(.7,2)(-.5,-.5)
\put(-3.9,0){\makebox(0,0){$\delta\left(\rule{0pt}{20pt}\right.$}}
\put(1.2,0){\makebox(0,0){$\left.\rule{0pt}{20pt}\right)$}}
\put(-.5,0){\makebox(0,0)[cc]{\Large$\nabla$}}
\put(-.5,.4){\vector(0,1){1}}
\put(.1,-1.25){\vector(-1,2){.5}}
\put(.67,-1.1){\vector(-1,1){1}}
\unitlength 3mm
\put(-21,-20){
\put(19.4,18){\vector(1,4){0.475}}
\put(18.9,18){\vector(1,2){.95}}
\put(16,18){\vector(2,1){3.8}}
\put(17.75,18){\makebox(0,0)[cc]{$\ldots$}}
\put(17.7,16.6){\makebox(0,0)[cc]{%
   $\underbrace{\rule{10mm}{0mm}}_{\mbox{\scriptsize $k$ inputs}}$}}
}
\end{picture}
\hskip 1cm := \hskip .5cm
G_k \hskip .2cm - \hskip -.2cm
\unitlength 3mm
\begin{picture}(20,4)(15.5,18.8)
\put(20,20.3){\vector(0,1){1.7}}
\put(18,18){\vector(1,1){1.8}}
\put(22,18){\vector(-1,1){1.8}}
\put(19,18){\vector(1,2){.88}}
\put(20,20){\makebox(0,0)[cc]{\Large$\circ$}}
\put(20.5,18){\makebox(0,0)[cc]{$\ldots$}}
\put(20,16.6){\makebox(0,0)[cc]{%
   $\underbrace{\rule{12 mm}{0mm}}_{\mbox{\scriptsize $k+2$}}$}}
\end{picture}
\raisebox{-2.5em}{\rule{0pt}{0pt}}
\end{equation}
where $G_k$ is a linear combination of $2$-vertex trees with one
$\nabla$-vertex~(\ref{nabla}), with $u < k$, and one white
vertex~(\ref{white}) with $u < k+2$. Explicit formulas for $k=0,1$ are
\[
\unitlength .5cm
\begin{picture}(.7,2)(-.5,-.5)
\put(-1.9,0){\makebox(0,0){$\delta\left(\rule{0pt}{20pt}\right.$}}
\put(1.2,0){\makebox(0,0){$\left.\rule{0pt}{20pt}\right)$}}
\put(-.5,0){\makebox(0,0)[cc]{\Large$\nabla$}}
\put(-.5,.4){\vector(0,1){1}}
\put(.1,-1.25){\vector(-1,2){.5}}
\put(.67,-1.1){\vector(-1,1){1}}
\end{picture} \hskip 1cm
:= \hskip .3cm -
\unitlength .3cm
\begin{picture}(4,4)(-2,-1)
\put(0,0.3){\vector(0,1){1.7}}
\put(-1,-2){\vector(1,2){.89}}
\put(1,-2){\vector(-1,2){.89}}
\put(0,0){\makebox(0,0)[cc]{\Large$\circ$}}
\end{picture}
\]
and
\[
\unitlength .5cm
\begin{picture}(.7,2)(-.5,-.5)
\put(-2.4,0){\makebox(0,0){$\delta\left(\rule{0pt}{20pt}\right.$}}
\put(1.2,0){\makebox(0,0){$\left.\rule{0pt}{20pt}\right)$}}
\put(-.5,0){\makebox(0,0)[cc]{\Large$\nabla$}}
\put(-.5,.4){\vector(0,1){1}}
\put(.1,-1.25){\vector(-1,2){.5}}
\put(.67,-1.1){\vector(-1,1){1}}
\put(-1.7,-1.1){\vector(1,1){1}}
\end{picture} \hskip 1cm
:= \hskip .3cm
\unitlength .3cm
\begin{picture}(4,4)(-2,-2)
\put(0,0.3){\vector(0,1){1.7}}
\put(-2,-4){\vector(1,2){1.8}}
\put(1.5,-1.5){\vector(-1,1){1.3}}
\put(0,0){\makebox(0,0)[cc]{\Large$\circ$}}
\end{picture}
\unitlength .5cm
\begin{picture}(.7,2)(-.2,0.1)
\put(-.5,0){\makebox(0,0)[cc]{\Large$\nabla$}}
\put(.1,-1.25){\vector(-1,2){.5}}
\put(.67,-1.1){\vector(-1,1){1}}
\end{picture} 
\hskip .5cm - \hskip 1cm
\begin{picture}(.7,2)(-.5,-1)
\put(-.5,0){\makebox(0,0)[cc]{\Large$\nabla$}}
\put(-.5,.4){\vector(0,1){1}}
\put(-.37,-.9){\vector(0,1){.7}}
\put(1.2,-2){\vector(-2,3){1.4}}
\end{picture}
\unitlength .3cm
\begin{picture}(4,4)(0.94,.3)
\put(-1.5,-1.5){\vector(1,1){1.3}}
\put(.9,-1.5){\vector(-1,2){.65}}
\put(0,0){\makebox(0,0)[cc]{\Large$\circ$}}
\end{picture}
\hskip -0.1cm - \hskip 1cm
\unitlength .5cm
\begin{picture}(.7,2)(-.5,-1)
\put(-.5,0){\makebox(0,0)[cc]{\Large$\nabla$}}
\put(-.5,.4){\vector(0,1){1}}
\put(.45,-1.9){\vector(-1,2){.83}}
\put(.58,-.8){\vector(-1,1){.83}}
\end{picture}
\unitlength .3cm
\begin{picture}(4,4)(-.9,-.05)
\put(-2.5,-1.5){\vector(3,2){2.22}}
\put(1.5,-1.5){\vector(-1,1){1.3}}
\put(0,0){\makebox(0,0)[cc]{\Large$\circ$}}
\end{picture}
\hskip -.1cm - \hskip .5cm
\unitlength .3cm
\begin{picture}(4,4)(-2,-1)
\put(0,0.3){\vector(0,1){1.7}}
\put(-2,-2){\vector(1,1){1.8}}
\put(2,-2){\vector(-1,1){1.8}}
\put(0,-2){\vector(0,1){1.75}}
\put(0,0){\makebox(0,0)[cc]{\Large$\circ$}}
\end{picture},
\raisebox{-1.8em}{\rule {0pt}{0pt}}
\]
which is a graphical form of an equation for the transformation of the
Christoffel symbols and their derivative under coordinate changes that can be
found in~\cite[Section~17.7]{kolar-michor-slovak} (but notice a different 
convention for covariant derivatives 
used in~\cite{kolar-michor-slovak}). 
\end{odstavec}

Finally, $\epsilon_v \in \{-1,+1\}$ in~(\ref{e68}) is a certain sign
whose definition can be found in~\cite[Section~4]{markl:na}. For the
purposes of this paper it will be enough to say that, if $\Graphname \in
\Gr^0(d)$ (no white vertices), then $\epsilon_v = 1$ for all $v \in
\Vert(\Graphname)$.  The relation between $\Gr^*(d)$ and natural
differential operators is described in:

\begin{theorem}[\cite{markl:na}]
\label{t2}
Each element in $H^0(\Gr^*(d),\delta) = \Ker\left(\delta : \Gr^0(d)\to
\Gr^1(d)\right)$ represents a~natural operator $\Con \times T^{\otimes d}
\to T$. On manifolds of dimension $\geq 2d-1$ this correspondence is
in fact an isomorphism
\[
H^0(\Gr^*(d),\delta) \cong \Nat(\Con \times T^{\otimes d},T).
\]
\end{theorem}

\begin{example}
In this example taken from~\cite{markl:na} we recall graphs
representing the Lie bracket, covariant derivative and
curvature. The Lie bracket  $[X,Y]$ of vector fields $X$, $Y$ 
is described~by
\[
\unitlength .7cm
\hskip -5cm \raisebox{.6cm}{$[X,Y]$:} \hskip 12mm
\begin{picture}(0,2.5)(0,0)
\put(0,2){\makebox(0,0)[cc]{\hskip .5mm$\bbox$}}
\put(0,1){\vector(0,1){.935}}
\put(0,1){\makebox(0,0)[cc]{\Large$\bullet$}}
\put(.4,1){\makebox(0,0)[lc]{\scriptsize$Y$}}
\put(0,0){\vector(0,1){.935}}
\put(0,0){\makebox(0,0)[cc]{\Large$\bullet$}}
\put(.4,0){\makebox(0,0)[lc]{\scriptsize$X$}}
\end{picture}
\hskip 1cm \raisebox{.6cm}{-}\hskip .8cm
\unitlength .7cm
\begin{picture}(0,2)(0,0)
\put(0,2){\makebox(0,0)[cc]{\hskip .5mm$\bbox$}}
\put(0,1){\vector(0,1){.935}}
\put(0,1){\makebox(0,0)[cc]{\Large$\bullet$}}
\put(.4,1){\makebox(0,0)[lc]{\scriptsize$X$}}
\put(0,0){\vector(0,1){.935}}
\put(0,0){\makebox(0,0)[cc]{\Large$\bullet$}}
\put(.4,0){\makebox(0,0)[lc]{\scriptsize$Y$}}
\end{picture}
\hskip 1cm \raisebox{.6cm}{$\in \Gr^0(2)$,}
\]
which in the obvious way expresses the local formula~(\ref{e3334}).
The covariant derivative is given by the graph
\[
\raisebox{-3em}{\rule{0pt}{0pt}}
\hskip 1.2cm
\unitlength .5cm
\begin{picture}(.7,3)(6,-.2)
\put(-.5,0){\makebox(0,0)[cc]{\Large$\nabla$}}
\put(2,0){\makebox(0,0)[c]{$+$}}
\put(-3,0){\makebox(0,0)[r]{$\nabla_X Y:$}}
\put(-.5,.4){\vector(0,1){1.1}}
\put(.1,-1.25){\vector(-1,2){.5}}
\put(.67,-1.1){\vector(-1,1){1}}
\put(-.45,1.4){\makebox(0,0)[bc]{$\bbox$}}
\put(0.2,-1.5){\makebox(0,0)[cc]{\Large$\bullet$}}
\put(0.2,-2){\makebox(0,0)[tc]{\scriptsize $X$}}
\put(1.5,-1.7){\makebox(0,0)[tc]{\scriptsize $Y$}}
\put(0.85,-1.3){\makebox(0,0)[cc]{\Large$\bullet$}}
\put(4,.3){
\unitlength .7cm
\begin{picture}(0,2.5)(0,1)
\put(0,2){\makebox(0,0)[cc]{\hskip .5mm$\bbox$}}
\put(0,1){\vector(0,1){.935}}
\put(0,1){\makebox(0,0)[cc]{\Large$\bullet$}}
\put(.4,1){\makebox(0,0)[lc]{\scriptsize$Y$}}
\put(0,0){\vector(0,1){.935}}
\put(0,0){\makebox(0,0)[cc]{\Large$\bullet$}}
\put(.4,0){\makebox(0,0)[lc]{\scriptsize$X$}}
\end{picture}
}
\end{picture}
\hskip -.3cm \raisebox{0.1cm}{$\in \Gr^0(2)$,}
\]
which is a graphical form of~(\ref{e10}).
Finally, the curvature $R: \Con \times T^{\otimes 3} \to T$ is given
by the graph
\[
\raisebox{-6em}{\rule{0pt}{0pt}}
\unitlength .5cm
\hskip 7cm
\begin{picture}(.7,3)(6,-.5)
\put(-.5,0){\makebox(0,0)[cc]{\Large$\nabla$}}
\put(2,0){\makebox(0,0)[c]{$-$}}
\put(-3,0){\makebox(0,0)[r]{$R(X,Y)(Z):$}}
\put(-.5,.4){\vector(0,1){1.1}}
\put(.1,-1.25){\vector(-1,2){.5}}
\put(.67,-1.1){\vector(-1,1){1}}
\put(-1.3,-1.12){\vector(2,3){.65}}
\put(-.45,1.4){\makebox(0,0)[bc]{$\bbox$}}
\put(0.2,-1.5){\makebox(0,0)[cc]{\Large$\bullet$}}
\put(0.2,-2){\makebox(0,0)[tc]{\scriptsize $X$}}
\put(-1.7,-1.9){\makebox(0,0)[tc]{\scriptsize $Y$}}
\put(1.5,-1.7){\makebox(0,0)[tc]{\scriptsize $Z$}}
\put(0.85,-1.3){\makebox(0,0)[cc]{\Large$\bullet$}}
\put(-1.4,-1.3){\makebox(0,0)[cc]{\Large$\bullet$}}
\end{picture}
\put(-.5,0){
\unitlength .5cm
\begin{picture}(.7,3)(1,-.5)
\put(-.5,0){\makebox(0,0)[cc]{\Large$\nabla$}}
\put(2,0){\makebox(0,0)[c]{$+$}}
\put(6.8,0){\makebox(0,0)[c]{$-$}}
\put(-.5,.4){\vector(0,1){1.1}}
\put(.1,-1.25){\vector(-1,2){.5}}
\put(.67,-1.1){\vector(-1,1){1}}
\put(-1.3,-1.12){\vector(2,3){.65}}
\put(-.45,1.4){\makebox(0,0)[bc]{$\bbox$}}
\put(0.2,-1.5){\makebox(0,0)[cc]{\Large$\bullet$}}
\put(0.2,-2){\makebox(0,0)[tc]{\scriptsize $Y$}}
\put(-1.7,-1.9){\makebox(0,0)[tc]{\scriptsize $X$}}
\put(1.5,-1.7){\makebox(0,0)[tc]{\scriptsize $Z$}}
\put(0.85,-1.3){\makebox(0,0)[cc]{\Large$\bullet$}}
\put(-1.4,-1.3){\makebox(0,0)[cc]{\Large$\bullet$}}
\put(5.5,-.5){
\begin{picture}(.7,3)(1,-.5)
\put(-.5,0){\makebox(0,0)[cc]{\Large$\nabla$}}
\put(-.5,.4){\vector(0,1){1.1}}
\put(.1,-1.25){\vector(-1,2){.5}}
\put(.67,-1.1){\vector(-1,1){1}}
\put(-.45,1.4){\makebox(0,0)[bc]{$\bbox$}}
\put(0.2,-1.5){\makebox(0,0)[cc]{\Large$\bullet$}}
\put(0.2,-2){\makebox(0,0)[tc]{\scriptsize $Y$}}
\put(.56,-1){\qbezier(0,0)(1,-1)(1,-2)}
\end{picture}
\put(1.1,-3.4){
\unitlength .5cm
\begin{picture}(.7,3)(1,-.5)
\put(-.5,0){\makebox(0,0)[cc]{\Large$\nabla$}}
\put(.1,-1.25){\vector(-1,2){.5}}
\put(.67,-1.1){\vector(-1,1){1}}
\put(0.2,-1.5){\makebox(0,0)[cc]{\Large$\bullet$}}
\put(0.2,-2){\makebox(0,0)[tc]{\scriptsize $X$}}
\put(1.5,-1.7){\makebox(0,0)[tc]{\scriptsize $Z$}}
\put(0.85,-1.3){\makebox(0,0)[cc]{\Large$\bullet$}}
\end{picture}
}}
\put(10.2,-.5){
\begin{picture}(.7,3)(1,-.5)
\put(-.5,0){\makebox(0,0)[cc]{\Large$\nabla$}}
\put(-.5,.4){\vector(0,1){1.1}}
\put(.1,-1.25){\vector(-1,2){.5}}
\put(.67,-1.1){\vector(-1,1){1}}
\put(-.45,1.4){\makebox(0,0)[bc]{$\bbox$}}
\put(0.2,-1.5){\makebox(0,0)[cc]{\Large$\bullet$}}
\put(0.2,-2){\makebox(0,0)[tc]{\scriptsize $X$}}
\put(.56,-1){\qbezier(0,0)(1,-1)(1,-2)}
\end{picture}
\put(1.1,-3.4){
\unitlength .5cm
\begin{picture}(.7,3)(1,-.5)
\put(-.5,0){\makebox(0,0)[cc]{\Large$\nabla$}}
\put(.1,-1.25){\vector(-1,2){.5}}
\put(.67,-1.1){\vector(-1,1){1}}
\put(0.2,-1.5){\makebox(0,0)[cc]{\Large$\bullet$}}
\put(0.2,-2){\makebox(0,0)[tc]{\scriptsize $Y$}}
\put(1.5,-1.7){\makebox(0,0)[tc]{\scriptsize $Z$}}
\put(0.85,-1.3){\makebox(0,0)[cc]{\Large$\bullet$}}
\end{picture}
}}
\end{picture}
}
\hskip 5.3cm \raisebox{0.1cm}{$\in \Gr^0(3)$.}
\]

\noindent 
We recommend as an exercise to verify that all the above graphs
belong to the kernel of $\delta$.
\end{example}

\section{Auxiliary results}
\label{s3}\label{sec5}
\label{auxiliary}
\setcounter{equation}{0}

The results of this purely technical section will be used, in
Section~\ref{cohomology}, to construct an explicit basis of the vector
space $H^0(\Gr^*(d),\delta)$ recalled in Section~\ref{rules}.  Let us
consider a~bicomplex $\bbB = (B^{*,*}, \delta= \deltah + \deltav)$,
with $B^{*,*} = \bigoplus_{p,q \in {\mathbb Z}} B^{p,q}$ and
differentials
\[
\deltah: B^{p,q}\to B^{p+1,q},\ \deltav: B^{p,q}\to B^{p,q+1}.
\] 
We require, as usual, that
\begin{equation}
\label{e6}
\deltav^2 = 0,\ \deltah^2 =0\ \mbox { and }\ \deltav\deltah + \deltah\deltav = 0.
\end{equation}
The associated total complex $\Tot(\bbB) = (B^*,\delta)$ has
$B^* := \bigoplus_{n \in {\mathbb Z}} B^{n}$ with $B^n :=
\bigoplus_{p+q=n} B^{p,q}$ and $\delta: = \deltah + \deltav : B^n \to
B^{n+1}$; see~\cite[\S XI.6]{maclane:homology} for the terminology.
Throughout this section we assume that
\begin{itemize}
\item[(i)] $\bbB$ is concentrated in the sector $0 \leq -p \leq q$ (see
Figure~\ref{f1}),

\item[(ii)] $\bbB$ is left-bounded in the sense 
that $B^{p,*} = 0$ for $p <\!\!< 0$, and\label{ii}

\item[(iii)]
the horizontal cohomology of $\bbB$ is concentrated on the
diagonal $p+q=0$, that is
\[
H^p(B^{*,q},\deltah) =0 \mbox{ for } p+q \not= 0
\]
or, equivalently, 
\[
H^m(B^*,\deltah) = 0 \mbox { for } m \not= 0.
\]
\end{itemize}

\begin{figure}[t]
\setlength{\unitlength}{5em}
\begin{center}
\begin{picture}(5,5.5)(-3,0)
\multiput(1,1)(0,1){4}{\makebox(0,0){$0$}}
\multiput(0,0)(-1,1){5}{\makebox(0,0){$0$}}
\def\pomoc#1#2{\put(#1,#2){\put(0,1){\makebox(0,0){$B^{#1,#2}$}}}}
\pomoc00
\pomoc01
\pomoc02
\pomoc03
\pomoc{-1}1
\pomoc{-1}2
\pomoc{-1}3
\pomoc{-2}2
\pomoc{-2}3
\pomoc{-3}3
\multiput(0,.25)(0,1){5}{\vector(0,1){.5}\put(.2,.25)%
          {{\makebox(0,0){$\deltav$}}}}
\multiput(-1,1.25)(0,1){4}{\vector(0,1){.5}\put(.2,.25)%
          {{\makebox(0,0){$\deltav$}}}}
\multiput(-2,2.25)(0,1){3}{\vector(0,1){.5}\put(.2,.25)%
          {{\makebox(0,0){$\deltav$}}}}
\multiput(-3,3.25)(0,1){2}{\vector(0,1){.5}\put(.2,.25)%
          {{\makebox(0,0){$\deltav$}}}}
\multiput(-4,4.25)(0,1){1}{\vector(0,1){.5}\put(.2,.25)%
          {{\makebox(0,0){$\deltav$}}}}
\multiput(-3.75,4)(1,0){5}{\vector(1,0){.5}\put(-.25,.2)%
          {{\makebox(0,0){$\deltah$}}}}
\multiput(-2.75,3)(1,0){4}{\vector(1,0){.5}\put(-.25,.2)%
          {{\makebox(0,0){$\deltah$}}}}
\multiput(-1.75,2)(1,0){3}{\vector(1,0){.5}\put(-.25,.2)%
          {{\makebox(0,0){$\deltah$}}}}
\multiput(-.75,1)(1,0){2}{\vector(1,0){.5}\put(-.25,.2)%
          {{\makebox(0,0){$\deltah$}}}}
\multiput(-3.75,4)(1,0){1}{\vector(1,0){.5}\put(-.25,.2)%
          {{\makebox(0,0){$\deltah$}}}}
\multiput(-4,5)(1,0){5}{$\vdots$}
\end{picture}
\caption{\label{f1}The relevant part of a bicomplex $\bbB$
  concentrated in the sector~\hbox{$0 \leq -p \leq q$.}}
\end{center}
\end{figure}

It then follows from a standard spectral sequence 
argument~\cite[Theorem~XI.6.1]{maclane:homology} that
\[
H^{m}(\Tot(\bbB)) = 0 \mbox{  for } m \not= 0 
\]
while $H^0(\Tot(\bbB))$ is isomorphic to the
direct sum $\Zh := \bigoplus_{r \geq 0}\Zh^r$ (which is finite,
by~(ii)) of subspaces
\begin{equation}
\label{Zh}
\Zh^r := \Ker(\deltah :B^{-r,r} \to B^{-r+1,r}). 
\end{equation}
Let us indicate  how to construct such an isomorphism.

\begin{proposition}
\label{p98}
Let $\beta : \Zh = \bigoplus_{r \geq 0}\Zh^r \to  \bigoplus_{r \geq 0}
B^{-r,r}$ be a linear map such that, for each $r \geq 0$ and $z \in
\Zh^r$, $\beta(z)$ is a cocycle in the total complex $\Tot(\bbB)$ and has
the form
\begin{equation}
\label{e56}
\beta(z) = z + \lot
\end{equation}
with some $\lot \in \bigoplus_{p > r} B^{-p,p}$. Then 
the map $\beta$ induces an isomorphism (denoted by the
same symbol)
\[
\beta: \Zh \stackrel{\cong}\to H^0(\Tot(\bbB)).
\]
\end{proposition}

\begin{proof}
Let us interpret $\Zh = \bigoplus_{r \geq 0}\Zh^r$ as a cochain 
complex concentrated in degree zero, with trivial differential.
Define a decreasing filtration 
\[
\Zh = F'_0 \Zh \supset  F'_1 \Zh \supset  F'_2 \Zh \supset  F'_3 \Zh
\supset \cdots 
\] 
of $\Zh$ by $F'_s \Zh :=
\bigoplus_{r \geq s} \Zh^r$. Similarly, define a decreasing filtration
\[
B^* = F''_0 B^* \supset  F''_1 B^* \supset  F''_2 B^* \supset  F''_3 B^*
\supset \cdots 
\]
of the total complex $\Tot(\bbB) = (B^*,d)$ by $F''_sB^* : =
\bigoplus_{p+q=*} \bigoplus_{q \geq s}B^{p,q}$. With these
definitions, $\beta$ is a map of filtered cochain complexes that
induces an isomorphism of the $E^1$-terms of the associated spectral
sequences. The proposition follows from a standard spectral sequence
argument~\cite[Theorem~XI.1.1]{maclane:homology}.
\end{proof}
 
In Proposition~\ref{p98}, l.o.t.~is an abbreviation from {\em lower
order terms\/}. The justification for this terminology will became
obvious on page~\pageref{gerda} of Section~\ref{cohomology}. 
It is not difficult to show that there
always exists a map $\beta$ satisfying the assumptions of the
proposition. For further applications we, however, need an
explicit construction of this map. 
It starts by choosing, for each $n \geq 0$, a
complementary subspace $D^n \subset B^{-n,n}$ to $\Zh^n$ so that
\begin{equation}
\label{e4}
B^{-n,n} = \Zh^n \oplus D^n.
\end{equation}
Let $\pi_n :  B^{-n,n} \to D^n$ be the projection. 
For $X^n := \Ker(\deltav\deltah : B^{-n,n} \to
B^{-n+1,n+1})$ we define $U : X^n \to X^{n+1}$ by
\begin{equation}
\label{e22}
U(x) := \pi_{n+1} \deltah^{-1} \deltav(x), \ x \in X^n.
\end{equation}

We must verify that this definition of $U$ makes sense. 
Since, for  $x \in X^n$, $\deltav \deltah x  = - \deltah \deltav x 
= 0$ (see the third equation of~(\ref{e6})), $\deltav x$ is a
$\deltah$-cocycle, so
$\deltav x \in \Im(\deltah :B^{-n-1,n+1} \to
B^{-n,n+1})$ by the acyclicity~(iii). The set $\deltah^{-1}
\deltav(x)$ is therefore non-empty. If $u',u'' \in \deltah^{-1}
\deltav(x)$, then $\deltah (u' - u'') = 0$ so
$\pi_{n+1} (u') =\pi_{n+1} (u'')$, which means that $U(x)=\pi_{n+1}
\deltah^{-1} \deltav(x)$ is a one-element set. 
The condition $\deltav\deltah U(x)=0$ follows from the simple fact that
\begin{equation}
\label{e1}
\deltah U(x) =  \deltah \pi_{n+1} \deltah^{-1} \deltav(x) = \deltav x.
\end{equation}
This shows that indeed $U(x) \in X^{n+1}$.
For $z \in \Zh^n \subset X^n$ finally define
\begin{equation}
\label{e2}
\beta (z) := z - U(z) + U^2(z) - U^3 (z) + \cdots = (\iden +U)^{-1}(z). 
\end{equation}
The above sum is, by assumption~(ii), finite.

\begin{lemma}
\label{l1}
For each $n \geq 0$ and $z \in \Zh^n$, the element 
$\beta(z) \in \bigoplus_{p
\geq n}B^{-p,p}$ is a degree $0$ cocycle of the total complex
$\Tot(\bbB)$. The map $\beta$ defined in~(\ref{e2}) therefore satisfies
assumptions of Proposition~\ref{p98}. 
\end{lemma}

\begin{proof}
By the definition of the differential in  
the total complex, one needs to verify that
\[
\deltah z =0,\
\deltah U(z) = \deltav z,\
\deltah U^2(z) = \deltav U(z), \ldots
\]
The above equations immediately follow from~(\ref{e1}).
\end{proof}

Let us formulate a simple lemma which will be used in the
proof of Corollary~\ref{c1} below.

\begin{lemma}
\label{l21}
For $x \in X^n$, the element $U(x)$ defined in~(\ref{e22})
is characterized by $U(x) \in D^{n+1}$ and \ $\deltah
U(x) = \deltav x$.
\end{lemma}

\begin{proof}
It follows from definition and
from~(\ref{e1}) that $U(x)$ satisfies the conditions of the 
lemma. On the
other hand,  suppose we are given an element $a \in D^{n+1}$ such that
$\deltah a = \deltav x$. Then $a = \pi_{n+1} a = \pi_{n+1} \deltah^{-1}\deltah a = 
\pi_{n+1} \deltah^{-1}\deltav x = U(x)$. 
\end{proof}

\begin{corollary}
\label{c1}
Let $G$ be a group and assume that the bicomplex $\bbB = (B^{*,*},\delta=
\deltah + \deltav)$ consists of reductive $G$-modules. Suppose moreover that
the differentials $\deltah$ and $\deltav$ are $G$-equivariant. Then there
exists a $G$-equivariant $\beta$ satisfying assumptions of
Proposition~\ref{p98}. 
\end{corollary}

\begin{proof}
By the reductivity of the actions, one may obviously assume that the
decomposition~(\ref{e4}) is $G$-invariant and that the projection
$\pi_n$ is $G$-equivariant. Moreover, $X^n$ is a $G$-stable subspace
of $B^{-n,n}$ and $\deltah (U(x)g) = (\deltah U(x))g = \deltav(x)g$
for each $x \in X^n$, $g \in G$, because the differentials are
assumed to be $G$-equivariant. Lemma~\ref{l21} then implies that
$U(xg) = U(x)g$, that is, $U$ is $G$-equivariant. Then $\beta$ defined
by~(\ref{e2}) is $G$-equivariant,~too.
\end{proof}

\section{Cohomology of graph complexes}
\label{cohomology}\label{sec6}
\setcounter{equation}{0}

In this part we apply the methods of Section~\ref{auxiliary} to the
graph complex constructed in Section~\ref{rules}. Let us start with
$\Gr^*(d)$. The first step is to realize that $(\Gr^*(d),\delta)$ is
the total complex of a bicomplex defined as follows.  For $p,q \in
{\mathbb Z}$, let
\begin{equation}
\label{bigrad}
\Gr^{p,q}(d) := \Span\left\{
\mbox {\rm graphs } \Graphname \in  \Gr^{p+q}(d);\
\mbox {\rm the number of } \nabla\mbox{\rm-vertices } = -p
\right\}.
\end{equation}
Define the horizontal differential $\deltah : \Gr^{p,q}(d) \to
\Gr^{p+1,q}(d)$ by
\begin{equation}
\label{e15}
\hskip 2.2cm
\deltah\left(\hskip 1.5cm
\unitlength .5cm
\begin{picture}(.7,2)(-.5,-.5)
\put(-.5,0){\makebox(0,0)[cc]{\Large$\nabla$}}
\put(-.5,.4){\vector(0,1){1}}
\put(.1,-1.25){\vector(-1,2){.5}}
\put(.67,-1.1){\vector(-1,1){1}}
\unitlength 3mm
\put(-21,-20){
\put(19.4,18){\vector(1,4){0.475}}
\put(18.9,18){\vector(1,2){.95}}
\put(16,18){\vector(2,1){3.8}}
\put(17.75,18){\makebox(0,0)[cc]{$\ldots$}}
\put(17.7,16.6){\makebox(0,0)[cc]{%
   $\underbrace{\rule{10mm}{0mm}}_{\mbox{\scriptsize $k$ inputs}}$}}
}
\end{picture}
\hskip .6cm
\right)
:= - \hskip -.3cm
\unitlength 3mm
\begin{picture}(20,4)(15.5,18.8)
\put(20,20.3){\vector(0,1){1.7}}
\put(18,18){\vector(1,1){1.8}}
\put(22,18){\vector(-1,1){1.8}}
\put(19,18){\vector(1,2){.88}}
\put(20,20){\makebox(0,0)[cc]{\Large$\circ$}}
\put(20.5,18){\makebox(0,0)[cc]{$\ldots$}}
\put(20,16.6){\makebox(0,0)[cc]{%
   $\underbrace{\rule{12 mm}{0mm}}_{\mbox{\scriptsize $k+2$}}$}}
\end{picture}
\raisebox{-1cm}{{\rule{0mm}{0mm}}}
\end{equation}
while $\deltah$ is trivial on remaining vertices. The vertical
differential $\deltav : \Gr^{p,q}(d) \to \Gr^{p,q+1}(d)$ is defined by
requiring that $\deltav : = \delta$ on black vertices~(\ref{black}), white
vertices~(\ref{white}) and the anchor~\anchor, while
\[
\deltav\left(\hskip 1.5cm
\unitlength .5cm
\begin{picture}(.7,2)(-.5,-.5)
\put(-.5,0){\makebox(0,0)[cc]{\Large$\nabla$}}
\put(-.5,.4){\vector(0,1){1}}
\put(.1,-1.25){\vector(-1,2){.5}}
\put(.67,-1.1){\vector(-1,1){1}}
\unitlength 3mm
\put(-21,-20){
\put(19.4,18){\vector(1,4){0.475}}
\put(18.9,18){\vector(1,2){.95}}
\put(16,18){\vector(2,1){3.8}}
\put(17.75,18){\makebox(0,0)[cc]{$\ldots$}}
\put(17.7,16.6){\makebox(0,0)[cc]{%
   $\underbrace{\rule{10mm}{0mm}}_{\mbox{\scriptsize $k$ inputs}}$}}
}
\end{picture}
\hskip .6cm
\right)
:= G_k, \hskip 3cm
\raisebox{-1.2cm}{{\rule{0mm}{0mm}}}
\]
where $G_k$ is the same as in~(\ref{e11}).

\begin{lemma}
\label{l22}
The object $\Gr^{*,*}(d) = (\Gr^{*,*}(d),\deltah + \deltav)$ constructed above
is a bicomplex
whose total complex is the graph complex $(\Gr^*(d),\delta)$
recalled in Section~\ref{rules}.
\end{lemma}

\begin{proof}
The only property which is not obvious are the relations~(\ref{e6})
which can be verified directly. 
\end{proof}

Let us check that the bicomplex $\Gr^{*,*}(d)$ satisfies conditions
(i)--(iii) on page~\pageref{ii}, 
Section~\ref{auxiliary}. One immediately sees that~(i)
is equivalent to the obvious inequality
\[
0 \leq
\mbox {\rm number of } \nabla \mbox { \rm vertices }
\leq 
\mbox {\rm number of } \nabla \mbox { \rm vertices } +
\mbox {\rm number of white vertices.}
\]
Simple graph combinatorics implies that 
each graph $\Graphname\in \Gr^*(d)$ has at most
$d-1$ $\nabla$-vertices, therefore  $\Gr^{p,q}(d) = 0$ for $p \leq
-d$ so the condition~(ii) of Section~\ref{auxiliary} is also satisfied.

To verify~(iii), we need to follow~\cite{markl:na} and observe that
$(\Gr^*(d),\deltav)$ is a particular case of 
the following construction.
{}For each collection $(E^*,\vt_E) = \{(E^*(s),\vt_E)\}_{s
\geq 2}$ of right dg-$\Sigma_s$-modules $(E^*(s),\vt_E)$, one
considers the complex
$\Grd *{E^*} = (\Grd *{E^*},\vt)$ spanned by graphs with 
$d$ black vertices~(\ref{black}), one vertex \anchor\ 
and a finite number of vertices decorated by elements
of $E$. The grading of $\Grd *{E^*}$ 
is induced by the grading of $E^*$ and the differential $\vt$ replaces
$E$-decorated vertices, one at a time, by their $\vt_E$-images, leaving other
vertices unchanged.
It is a standard fact~\cite{mv} (see also~\cite[Theorem~21]{markl:ha}) 
that the assignment 
$(E^*,\vt_E) \mapsto (\Grd *{E^*},\vt)$ 
is a polynomial, hence exact, functor, thus
\begin{equation}
\label{21}
H^*(\Grd *{E^*},\vt) \cong \Grd *{H^*(E,\vt_E)}.
\end{equation}

Let now $(E^*,\vt_E) = \{(E^*(s),\vt_E)\}_{s \geq 2}$ be such that
$E^0(s)$ is spanned by symbols~(\ref{nabla}) with $u+2 = s$, $E^1(s)$
by symbols~(\ref{white}) with $u = s$, and $E^m(s) = 0$ for $m \geq
2$. The differential $\vt_E$ is defined by replacement
rule~(\ref{e15}).  More formally, 
\[
E^0(s) = {\rm Ind}^{\Sigma_s}_{\Sigma_{s-2}\times \Sigma_2} 
({{\mathbf 1}}_{s-2}\times {{\mathbf 1}}_2)
\mbox { and } E^1(s) = {{\mathbf 1}}_s,
\] 
where ${{\mathbf 1}}_{s-2}$ (resp.~${{\mathbf 1}}_{2}$,
resp.~${{\mathbf 1}}_s$) denotes the trivial representation of the
symmetric group $\Sigma_{s-2}$ (resp.~$\Sigma_2$,
resp.~$\Sigma_s$). The differential $\vt_E$ then sends the generator
$1 \times 1\in {{\mathbf 1}}_{s-2}\times {{\mathbf 1}}_2$ into $-1 \in
{{\mathbf 1}}_{s}$.  It is clear that, with this particular choice of
the collection $(E^*,\vt_E)$,
\[
(\Gr^*(d),\deltah) \cong (\Grd *{E^*},\vt).
\]

Since $\vt_E : E^0(s) \to E^1(s)$ is onto, the collection
$H^*(E,\vt_E) = \{H^*(E(s),\vt_E)\}_{s \geq 2}$ is concentrated in degree 
$0$, with $H^0(E(s),\vt_E)$ the kernel
\begin{equation}
\label{kernel}
\Kr(s) : = \Ker\left(\vt_E : E^0(s) \to E^1(s)\right).
\end{equation}
\label{kriz}
Denoting by $\Kr$ the collection $\Kr := \{\Kr(s)\}_{s \geq 2}$ we
conclude that 
\begin{equation}
\label{Jark}
H^*(\Gr^*(d),\deltah) \cong \Gr[\Kr](d).
\end{equation}
The right hand side is concentrated in degree zero so
we omitted the star indicating the grading. In particular,
$H^m(\Gr^*(d),\deltah) = 0$ for $m \not= 0$ which establishes~(iii).

The above calculation shows that, for the bicomplex $\Gr^{*,*}(d)$,
the cocycle space $\Zh^n$ of~(\ref{Zh}) equals
\begin{equation}
\label{Zh1}
\Zh^n= \Span\left\{
\Graphname \in \Gr[\Kr](d);\ \mbox {$\Graphname$ has precisely $n$ vertices
  decorated by $\Kr$}\right\}.
\end{equation}

Let $\alpha_\sigma$, $\sigma \in {\rm Ush}(s-2,2)$, be real
coefficients as in Theorem~A.
If we take the symbol~(\ref{nabla}), with the inputs numbered
consecutively from left to right by $\{\rada 1s\}$, as the generator of
$E^0(s)$, then $\Kr(s)$ is, as a $\Sigma_s$-module, generated by the
linear combination
\begin{equation}
\label{xi}
\xi_s :=  \sum_{\sigma \in {\rm Ush}(s-2,2)} \hskip 4.2em
\unitlength 5mm
\begin{picture}(.7,2)(-.5,-.5)
\put(-.5,0){\makebox(0,0)[cc]{\Large$\nabla$}}
\put(-.5,.4){\vector(0,1){1}}
\put(.1,-1.25){\vector(-1,2){.5}}
\put(.59,-1.22){\vector(-1,1){.97}}
\unitlength 3mm
\put(-21,-20){
\put(19.4,18){\vector(1,4){0.475}}
\put(18.9,18){\vector(1,2){.95}}
\put(16,18){\vector(2,1){3.8}}
\put(18.75,17){\makebox(0,0)[cc]{$\ldots$}}
}
\put(-6.5,-3.2){{\scriptsize $\sigma(1)$}}
\put(-.1,-3.2){{\scriptsize $\sigma(n)$}}
\end{picture}
\hskip 1.5em =   \sum_{\sigma \in {\rm Ush}(s-2,2)} \hskip 3.2em
\begin{picture}(.7,2)(-.5,-.5)
\put(-.5,0){\makebox(0,0)[cc]{\Large$\nabla$}}
\put(-.5,.4){\vector(0,1){1}}
\put(.1,-1.25){\vector(-1,2){.5}}
\put(.59,-1.22){\vector(-1,1){.97}}
\unitlength 3mm
\put(-21,-20){
\put(19.4,18){\vector(1,4){0.475}}
\put(18.9,18){\vector(1,2){.95}}
\put(16,18){\vector(2,1){3.8}}
\put(17.75,18){\makebox(0,0)[cc]{$\ldots$}}
}
\end{picture}
\hskip .3em
\cdot \sigma\ .
\raisebox{-2.5em}{\rule{0pt}{0pt}}
\end{equation}

For each $n \geq 0$ consider the subcomplex $\Gr^*(n+1)_n$ of
$\Gr^*(n+1)$ spanned by graphs with a distinguished subset  
$\{\anteanchor1, \ldots, \anteanchor n \hskip .1em \}$ 
of the set of black vertices~(\ref{black}) with $u=0$,
see~\ref{ss2}. Suppose we are given, for each $n \geq 0$, cochains
$\nu_n \in \Gr^0(n+1)_n$ of the form
\begin{equation}
\label{nn}
\nu_n = b_n + \lot
\end{equation}
where $b_n$ denotes the black vertex~(\ref{black}) with\label{gerda}
$u=n$. The abbreviation l.o.t.\ denotes a linear combination of graphs in
$\Gr^0(n+1)_n$ that has at least one $\nabla$-vertex. It is not
difficult to verify that each such a graph represents a local
expression whose differential order is strictly smaller than the
differential order of the local expression represented by $b_n$ (which
is $n-1$). This explains why l.o.t.\ abbreviates ``lower order terms.''

Similarly, recall from~\ref{ss2} that $\Gr^*(n)_n$ denotes the
subcomplex of $\Gr^*(n)$ spanned by graphs whose all black vertices
belongs to the distinguished subset  
$\{\anteanchor1, \ldots, \anteanchor n \hskip .1em\}$.
Suppose that we are given, for each $n \geq 3$, cochains $\varsigma_n \in
\Gr^0(n)$ of the form
\begin{equation}
\label{kn}
\varsigma_n = \xi_n + \lot
\end{equation}
where $\xi_n$ is as in~(\ref{xi}) and l.o.t.\ a linear
combination of graphs with at least two $\nabla$-vertices. It is not
difficult to see that $\nu_n$'s and $\varsigma_n$'s as above always
exist, but we prove a stronger result:

\begin{proposition}
\label{eqv}
The cocycles $\{\nu_n\}_{n \geq 2}$ and $\{\varsigma_n\}_{n \geq 3}$ can
be choosen `equivariantly,' that is, in such a way that they enjoy the
same symmetries as the elements $\{b_n\}_{n \geq 2}$  and
$\{\xi_n\}_{n \geq 3}$.
\end{proposition}

\begin{proof}
Repeating the reasonings in the proof of Lemma~\ref{l22}, one easily sees 
that the obvious modification of the bigrading~(\ref{bigrad}) turns
the graph complex $\Gr^*(n+1)_n$ into a bicomplex satisfying conditions
(i)--(iii) on page~\pageref{ii}. The symmetric group
$\Sigma_n$ permutes the distinguished vertices of graphs spanning
$\Gr^*(n+1)_n$. This action satisfies the requirements of
Corollary~\ref{c1} which therefore gives a $\Sigma_n$-equivariant 
$\beta$ satisfying assumptions of
Proposition~\ref{p98}. The element $\nu_n := \beta(b_n)$ is then
an `equivariant' cocycle in that it is, as $b_n$,
$\Sigma_n$-stable. An `equivariant' $\varsigma_n$ can be constructed in
the same fashion, considering $\Gr^*(n)_n$ instead of $\Gr^*(n+1)_n$.
\end{proof}

\section{Proofs of Theorems~A--F}
\label{proofs}\label{sec7}
\setcounter{equation}{0}

\noindent 
{\bf Proof of Theorem~A.}  
We will in fact give two proofs, one using the classical reduction
theorem, and one based on the graph complex method. Let us start with
the `classical' proof.

By Example~\ref{e88}, expression~(\ref{konec}) corresponding to
the leading term of the\label{J}
iterated covariant derivative $\nabla^{n-2} R$ of the curvature tensor
belongs to the kernel $\Kr(n)$, for each $n \geq 3$. Since, by
assumption, expressions~(\ref{konec}) corresponding to 
the leading terms of the operators $D_n$ generate
$\Kr(n)$, one clearly has, for each $n \geq 3$,
\[
\nabla^{n-2} R(X_1,\dots,X_n) = \sum_{\sigma\in \Sigma_n}  A_\sigma 
D_{n}(X_{\sigma(1)},\dots,X_{\sigma(n)}) + \lot 
\]
with some real coefficients $A_\sigma$. Similarly, the leading term  of
any operator $V_n$ as in Theorem~A equals, 
by Example~\ref{e88} and~(\ref{Eq2.4}), 
the leading term of the iterated
covariant derivative $(\nabla^{n-1} X_n)(X_1,\dots,X_{n-1})$, therefore
\[
(\nabla^{n-1} X_n)(X_1,\dots,X_{n-1}) = V(\Gamma,X_1,\dots,X_n) + \lot
\] 
Theorem~A now follows from
Theorem~\ref{Th1.2} and obvious induction on the differential degree.

Let us give another proof of Theorem~A based on the method of graph
complexes. The drawback of this proof is that it requires the
`stability' $\dim(M) \geq 2d-1$. On the other hand, it is completely
independent on local calculations. Moreover, we will need to set up
the stage for graph-ical proofs of other statements.

Each iteration in the sense of~\ref{iter}
is clearly a linear combination of terms given by contracting `free'
indices of the local coordinate expressions of the operators
$\{D_n\}_{n \geq 3}$ and $\{V_n\}_{n \geq 2}$. Each such a contraction
is determined by a `contraction scheme,' which is a graph as
in~\ref{o1} with vertices of the following two types:

\begin{itemize}
\item[--] 
vertices $d_n$, $n \geq 3$, with $n$ linearly
ordered input edges and one output, and
\item[--] 
vertices $v_n$, $n \geq 0$, labeled $\rada 1d$, with $n$ linearly
ordered edges and one output.
\end{itemize}

Observe that we allowed vertices $v_n$ also for $n=0$ as places where order
zero vector field variables are to be inserted.
Denote by $\Cont(d)$ the space spanned by the above contraction
schemes. Consider the diagram
\begin{equation}
\label{diag}
\Gr[\Kr](d) \stackrel\pi\twoheadleftarrow \Cont(d)  \stackrel\Psi\to 
\Gr^0(d)
\end{equation}
in which the maps $\pi$ and $\Psi$ are defined as follows. 

The map $\pi$ replaces each vertex $d_n$ of a contraction scheme $K
\in \Cont(d)$ by $\xi_n$ and each vertex $v_n$ by $b_n$ --
recall~(\ref{xi}) resp.~(\ref{black}) for the definition of $\xi_n$
resp.~$b_n$.  The map $\Psi$ replaces each vertex $d_n$ by the cocycle
$\varsigma_n \in \Gr^0(n)$ representing, in the correspondence of
Theorem~\ref{t2}, the operator $D_n$, and each vertex $v_n$ by the
cocycle $\nu_n \in \Gr^0(n+1)$ representing the operator $V_{n+1}$.
Therefore $\Psi(K)$ is the cocycle representing the iteration
determined by $K$.

The map $\pi$ is an epimorphism. One may establish this fact by
constructing a section (=~right inverse) $s : \Gr[\Kr](d) \to
\Cont(d)$ of $\pi$ as follows. Recall~\cite{markl:handbook} that a
graph with vertices decorated by a collection 
$F = \{F(s)\}_{s \geq 0}$ of right $\Sigma_s$-modules $F(s)$
($F$-graph for short) is an equivalence
class of graphs whose vertices have linearly ordered inputs and are
decorated by elements of $F$.  The equivalence
identifies graphs that differ only by the orders of the inputs and
actions of the corresponding permutations at the decorations. The
space $\Gr[\Kr](d)$ is then spanned by $F$-graphs with $F := \Kr
\oplus 1$, where $\Kr$ defined in~(\ref{kernel}) is generated by the
sequence $\{\xi_n\}_{n \geq 3}$ (we of
course put $\Kr(s) = 0$ for $s = 0,1,2$) and $1$ is the
collection of trivial representations generated by the elements
$\{b_n\}_{n \geq 0}$.

Assume that the graphs $\rada {\Graphname_1}{\Graphname_b}$ form a basis of
$\Gr[\Kr](d)$. Choose a representative $\widetilde{\Graphname}_i$ of each
$\Graphname_i$, $1 \leq i \leq b$, in the equivalence relation described
in the previous paragraph. Define $s (\Graphname_i)$ as the contraction
scheme obtained from $\widetilde{\Graphname}_i$ by replacing each
vertex $\xi_n$ by $d_n$ and each vertex $b_n$ by $v_n$, preserving
the linear orders of the inputs.
The identity $\Psi \circ s = {\it id}$ for the map $s$
defined in this way is obvious.

The composition $\beta := \Psi \circ s : \Gr[\Kr](d) \to \Gr^0(d)$ is
easy to describe; \label{pg2} $\beta(\Graphname_i)$ is the graph obtained
from ${\widetilde \Graphname}_i$ by replacing each vertex $\xi_n$ by the
graph $\varsigma_n$ representing the operator $D_n$ and each vertex $b_n$
of ${\Graphname}_i$ by the graph $\nu_n$ representing the operator $V_{n+1}$,
$1 \leq i \leq b$.  One easily sees that $\beta$ satisfies assumptions
of Propositions~\ref{p98}, with $B^{*,*}$ the graph bicomplex
$\Gr^{*,*}(d)$ defined by~(\ref{bigrad}) and $\Zh$ equalling,
by~(\ref{Zh1}), the space $\Gr[\Kr](d)$.  Therefore $\beta = \Psi
\circ s$ induces an isomorphism $\Gr[\Kr](d) \cong
H^0(\Gr^*(d),\delta)$. In particular, the map $\Psi$ is an epimorphism
onto $\Ker(\delta : \Gr^0(d) \to \Gr^1(d)) = H^0(\Gr^*(d),\delta)$.
This, along with Theorem~\ref{t2}, proves Theorem~A.

\vskip .2em

\noindent 
{\bf Proof of Theorem~B.}  Operators $\{D_n\}_{n \geq 3}$ defined
by~(\ref{j12}) are clearly iterations, in the sense of~\ref{iter}, of
the operators $\{U_n\}_{n \geq 3}$. If the leading terms of the
operators $\{D_n\}_{n \geq 3}$ are as in
Theorem~A, then each operator $\gO$ is an iteration of the operators
$\{D_n\}_{n \geq 3}$ and hence also of the operators  
$\{U_n\}_{n \geq 3}$. This proves (ii) $\Longrightarrow$ (i). 
Let us prove the oposite implication.

Fix $n \geq 3$ and write $U_n = L_n + \lot$ (i.e., $L_n$ is the
leading term of $U_n$).  Let $E_n = Q_n + \lot$ be an arbitrary
operator whose leading term satisfies the assumptions of Theorem~A. If~(ii) is
fulfilled, then, in particular, the operator $E_n$ is an iteration of
the operators $\{U_u\}_{u \leq n}$. A simple reasoning based on the
c-order implies that the leading term $Q_n$ of $E_n$ is obtained from
$L_n$ by linear combinations of successive applications of the
operations $\gO \mapsto \gO\sigma$, $\sigma \in \Sigma_n$, and $\gO
\mapsto \Tr_j(\gO) X_j$, with some $1 \leq j \leq n$, i.e.~by using
elementary iterations~(i), (ii) and (iv) of~\ref{iter} only, with
$\gO' := \gO$ and $\gO''(\Gamma,\Rada X1n) := X_j$ in (iv).

Traces commute with the symmetric group action in
the sense that, for each $\sigma \in \Sigma_n$ and $1 \leq j \leq
n$, there exists some
$\tilde \sigma_j \in \Sigma_n$ such that
\begin{equation}
\label{j19}
(\Tr_j(\gO) X_j) \sigma = \Tr_{\sigma^{-1}(j)}(\gO \tilde \sigma_j)
X_{\sigma^{-1}(j)}.
\end{equation} 
The explicit description of the permutation
$\tilde \sigma_j$ is not important for this proof and we leave it 
as an exercise for the reader. The `commutativity'~(\ref{j19})
implies that one may always move the symmetric group action inside the
trace and write
\begin{equation}
\label{e97}
Q_n = L_n c + \sum_{1 \leq j \leq n}\Tr_j(L_n c_j) X_j
\end{equation}
for some $c,\Rada c1n \in \RRR [\Sigma_n]$.

While the leading terms of natural differential operators need not be
invariant under general coordinate changes, they are still invariant
under the action of the general linear group ${\rm GL}_n$. This means
(see~\cite{markl:ig}) that the expressions in~(\ref{e97}) are
represented by linear combination of graphs.  These graphs have the
anchor \anchor, vertices $\{\anteanchor1, \ldots, \anteanchor n\}$,
and one vertex~(\ref{nabla}) with $u := n-2$.  Each such a graph has
$n+1$ edges, but, since the $\nabla$-vertex is fully symmetric in the
first $(n-2)$ and the last two inputs, the stability dimension
($=$~the minimal dimension of the underlying space for which a ${\rm
GL}_n$-invariant 
operator uniquely determines a linear combination of graphs) is $3$,
see~\cite[Proposition~4.9]{markl:ig}.

So both sides of~(\ref{e97}) are represented by the same linear
combinations of graphs, i.e.~they are given
by the same contraction schemes of indices. 
We conclude that, if $\dim(M) \geq 3$, the leading term of the
operator $D_n$ defined by~(\ref{j12}), with $c,\Rada c1n$ as
in~(\ref{e97}), is of
the form required by Theorem~A.

\vskip .2em

\noindent 
{\bf Proof of Theorem~C.}  One may assign to each graph $\Graphname \in
\Gr[\Kr](d)$ the (formal) vf-order (where, as on page~\pageref{huhu},
vf abbreviates `vector field') defined by the summation
\begin{equation}
\label{eqaut}
\vforder(\Graphname) := \sum_{v \in \Vert(\Graphname)} \vforder(v),
\end{equation}
where 
\[
\vforder(v) := \cases{0}{if $v$ is $\xi_n$, $n \geq 3$, and}
                        n{if $v$ is $b_n$, $n \geq 0$.}
\]

The vf-order of a contraction scheme $G \in \Cont(d)$ can be
defined similarly, with the role of vertices $\xi_n$ played by
$d_n$, and the role of vertices $b_n$ by $v_n$.
Therefore, if a contraction scheme has vertices
$\rada{v_{p_1}}{v_{p_t}}$ for some $\Rada
p1t \geq 0$ (plus possibly some other vertices of either types), then
\begin{equation}
\label{e45}
p_1 + \cdots + p_t \leq \vforder(G).
\end{equation}

Finally, the vf-order of a graph $\Graphname$ in $\Gr^0(d)$ is
given by formula~(\ref{eqaut}) in which we define now
\[
\vforder(v) := \cases{0}{if $v$ is a $\nabla$-vertex, and}
                        n{if $v$ is $b_n$, $n \geq 0$.}
\]

The vf-order of an element of $\Gr[\Kr](d)$ (resp.~$\Cont(d)$,
resp.~$\Gr^0(d)$) is then the maximum of
vf-orders of its linear constituents. 
It is clear that the (formal) vf-order of a cocycle in $\Gr^0(d)$
equals the vf-order of the operator it represents.

We are going to show that the isomorphism $\beta = \Phi \circ s$ constructed in
the proof of Theorem~A preserves the vf-order. As before, let
$\{\Graphname_i\}_{1 \leq i \leq b}$ be a basis of $\Gr[\Kr](d)$. Recall
from page~\pageref{pg2} that $\beta$ acts by replacing
$\xi_n$-vertices of $\Graphname_i$ by $\varsigma_n$ and $b_n$-vertices of
$\Graphname_i$ by $\nu_n$.
Another observation we need is that
\begin{equation}
\label{j1}
\vforder(\varsigma_n) = 0,\ n \geq 3,
\end{equation}
while $\nu_n = b_n + \eta_n$, where $\eta_n$ is a graph such that 
\begin{equation}
\label{j2}
\vforder(\eta_n) < \vforder(b_n),\ n \geq 0.
\end{equation}
Equation~(\ref{j1}) expresses that $\varsigma_n$
represents the operator $D_n$ which is, by assumption, of
order $0$ in the vector field variables. Inequality~(\ref{j2}) is a
consequence of the fact that $\eta_n$ is a linear combination of
graphs having at least one $\nabla$-vertex and that, by
simple graph combinatorics, the vf-order of each such a graph
is strictly less than the vf-order of $b_n$, compare the remark
following~(\ref{nn}). 

The above implies that $\beta(\Graphname_i) = \Graphname_i + \epsilon_i$, where
$\vforder(\epsilon_i) < \vforder(\Graphname_i)$ for each $1 \leq i \leq
b$, thus also, for an arbitrary linear combination 
$y \in \Gr[\Kr](d)$ of the basis elements, $\beta(y) = y +
\varepsilon$, where $\vforder(\varepsilon) < \vforder(y)$. We conclude
that then indeed $\vforder(y) = \vforder(\beta(y))$. The fact that the
section $s$ constructed in the proof of Theorem~A also preserves the
vf-order, is obvious.

Let $\gO \in \Nat(\Con \times T^{\otimes d},T)$ be a
differential operator represented by a cocycle $c \in \Gr^0(d)$, $y :=
\beta^{-1}(c)$ and $C := s(y)$. Acording to our constructions, $C \in
\Cont(d)$ describes an iteration of $\{D_n\}_{n \geq 3}$ and
$\{V_n\}_{n \geq 1}$ representing $\gO$.  Since both $\beta$ and $s$
preserve the vf-order, one has $\vforder(C) = \vforder(\gO)$. 
Theorem~C now immediately follows from formula~(\ref{e45}).

\vskip .2 em

\noindent 
{\bf Remaining proofs.}  A moment's reflection convinces us that the
kernel of the map $\Psi$ in diagram~(\ref{diag}) is generated by
contraction schemes expressing relations~(\ref{dev}) and (\ref{b5}).
This is precisely the content of Theorem~D.
The `ideal' tensors in Theorem~E are natural operators corresponding
to the cocycles $\{\varsigma_n\}_{n \geq 3}$ and
$\{\nu_n\}_{n \geq 2}$ constructed in Proposition~\ref{eqv}. Theorem~F
is a combination of Theorem~\ref{t2} and isomorphism~(\ref{Jark}).


\def\cprime{$'$}

\end{document}